\documentclass[journal ]{new-aiaa}

\let\transp\relax

\def\nobiblatex{}
\def\nohref{}
\def\noindentadjust{}

\usepackage{ifthen}
\usepackage{pgffor}
\usepackage{etoolbox}
\usepackage{parselines}
\usepackage{environ}
\usepackage{keyval}
\usepackage{listofitems}
\usepackage{forloop}
\usepackage{calc}
\usepackage{xstring}
\usepackage{changepage}
\usepackage{scrextend}
\usepackage{expl3}
\usepackage{stackengine}

\makeatletter
\ifthenelse{\isundefined{\lefteqs}}{
\usepackage{amsmath}
\usepackage{amsfonts}
\usepackage{mathtools}
}{
\usepackage[fleqn]{amsmath}
\usepackage[fleqn]{mathtools}
}
\usepackage{amssymb}
\@ifclassloaded{autart}{}{\usepackage{amsthm}}
\makeatother
\usepackage{thmtools}
\usepackage{ar}
\usepackage{empheq}
\usepackage{mathrsfs}
\usepackage{bm}
\usepackage{sansmath}

\makeatletter
\let\ltxxlabel\ltx@label
\makeatother

\usepackage{color}
\makeatletter
\@ifclassloaded{beamer}{}{
\usepackage[table]{xcolor}}
\makeatother

\usepackage{graphicx}
\makeatletter
\@ifclassloaded{uwthesis/uwthesis}{
  \usepackage{subfig}
}{
  \@ifclassloaded{IEEEtran}{
  }{
    \usepackage{subcaption}
  }
}
\@ifclassloaded{IEEEtran}{}{\usepackage{caption}}
\makeatother
\usepackage{dblfloatfix}

\usepackage{tabularx}
\usepackage{adjustbox}
\usepackage{multicol}
\usepackage{diagbox} 
\usepackage{varwidth} 
\usepackage{makecell} 
\usepackage{arydshln}

\usepackage{tikz}
\usepackage{tikz-3dplot}
\usetikzlibrary{backgrounds}
\usetikzlibrary{arrows}
\usetikzlibrary{arrows.meta}
\usetikzlibrary{intersections}
\usetikzlibrary{math}
\usetikzlibrary{fpu}
\usetikzlibrary{shapes}
\usetikzlibrary{positioning}
\usetikzlibrary{decorations.markings}
\usetikzlibrary{decorations.pathreplacing}
\usetikzlibrary{decorations.text}
\usetikzlibrary{tikzmark}
\usetikzlibrary{fadings}

\makeatletter
\@ifclassloaded{standalone}{}{\usepackage[tikz]{mdframed}}
\@ifclassloaded{beamer}{}{
  \@ifclassloaded{uwthesis/uwthesis}{}{
    \@ifclassloaded{IEEEtran}{}{
      \usepackage{tocloft}
    }}}
\usepackage{enumitem}
\makeatother
\usepackage{listings}
\usepackage{lettrine}
\usepackage{algorithm, algorithmicx}
\ifthenelse{\isundefined{\algouseends}}{
\usepackage[noend]{algpseudocode}}{
\usepackage{algpseudocode}}
\usepackage[many,magazine]{tcolorbox}
\usepackage[binary-units=true]{siunitx}
\usepackage{chemformula}
\usepackage[normalem]{ulem}
\usepackage{cancel}
\usepackage{accents}
\makeatletter
\@ifclassloaded{beamer}{\usepackage{multimedia}}{}
\makeatother
\usepackage{lipsum}
\usepackage{mfirstuc}
\usepackage{import}
\usepackage{nicematrix}

\ifthenelse{\isundefined{\nobiblatex}}{
  \usepackage{bibentry} 
  \newcommand{\usebiblatex}[1]{%
    \usepackage[backend=biber,
    style=#1,
    sorting=none,
    url=false,
    sortcites=true,
    defernumbers=true]{biblatex}}
  \ifthenelse{\isundefined{\preprintcsm}}{%
    \usebiblatex{numeric}}{%
    \usebiblatex{ieee}}
}{}

\usepackage[outline]{contour}
\usepackage{extdash}
\usepackage{textcomp}

\usepackage{float}
\usepackage{afterpage}
\ifthenelse{\isundefined{\landscapemode}}{
\def\landscapemode{pdflscape}
}{}
\usepackage{\landscapemode} 
\usepackage{varwidth}
\usepackage{nomencl}
\usepackage{longtable}
\usepackage{marginnote}
\makeatletter
\@ifclassloaded{beamer}{}{
\@ifclassloaded{autart}{}{
\@ifclassloaded{standalone}{}{
\@ifclassloaded{ieeeconf}{}{
\usepackage[explicit]{titlesec}
}}}}
\makeatother
\usepackage{soul}
\usepackage{stackengine}
\usepackage{scalerel}
\usepackage{lastpage}
\usepackage{fancyhdr}
\usepackage{nameref}

\makeatletter

\ifthenelse{\isundefined{\preprintcsm}}{}{
  \gdef\nohref{}}
\@ifclassloaded{beamer}{%
  \gdef\nohref{}}{}

\ifthenelse{\isundefined{\nohref}}{
  \usepackage[hypertexnames=false]{hyperref}}{}

\makeatother

\usepackage{xr}

\usepackage{xspace}
\usepackage{setspace}

\newcommand{\myappend}[1]{\ifthenelse{\equal{#1}{}}{}{,#1}}
\newcommand{\myprepend}[1]{\ifthenelse{\equal{#1}{}}{}{#1,}}

\ExplSyntaxOn
\makeatletter
\NewDocumentCommand{\parengen}{m m m m}
{
  \group_begin:
  \tl_set:Nn \l_left_paren_tl {#1}
  \tl_set:Nn \l_right_paren_tl {#2}
  \tl_set:Nn \l_paren_sz_tl {#3}
  \ifthenelse{\equal{\l_paren_sz_tl}{}}{
    \l_left_paren_tl #4 \l_right_paren_tl
  }{
    \ifthenelse{\equal{\l_paren_sz_tl}{auto}}{
      \left\l_left_paren_tl #4 \right\l_right_paren_tl
    }{
      \csname \l_paren_sz_tl\endcsname\l_left_paren_tl
      #4
      \csname \l_paren_sz_tl\endcsname\l_right_paren_tl
    }
  }
  \group_end:
}

\newcommand{\pare}[2][]{\parengen{(}{)}{#1}{#2}}
\newcommand{\brak}[2][]{\parengen{[}{]}{#1}{#2}}

\newcommand{\abso}[2][]{\parengen{|}{|}{#1}{#2}}

\newcommand{\@@normgen}[2][]{\parengen{\|}{\|}{#1}{#2}}

\makeatother
\ExplSyntaxOff

\let\oldbig\big

\let\oldBigg\Bigg
\renewcommand{\big}{\oldbig}

\renewcommand{\Bigg}{\oldBigg}

\makeatletter
\newcommand{\fade}[2][1]{%
  \edef\fade@amount{\the\numexpr 100-#1*30 \relax}%
  {\color{black!\fade@amount!white}#2}}
\makeatother

\newcommand{\ifelsemath}[3]{%
  \pgfmathparse{#1 ? 1:0}%
  \ifthenelse{\pgfmathresult>0}{#2}{#3}}

\newcommand{\inv}[1][]{\ifthenelse{\equal{#1}{}}{^{-1}}{^{-#1}}}

\ExplSyntaxOn
\makeatletter
\NewDocumentCommand{\norm}{o m}
{
  \group_begin:
  \IfNoValueTF{#1}{
    \def\@nargin{0}
  }{
    \seq_set_split:Nnn \l_var_seq { , } { #1 }
    \edef\@nargin{\seq_count:N \l_var_seq}
  }
  \seq_pop_left:NN \l_var_seq \l_norm_type_tl
  \seq_pop_left:NN \l_var_seq \l_delimiter_size_tl
  \pgfmathparse{\@nargin<=1 ? 1:0}
  \ifthenelse{\pgfmathresult>0}{
    \tl_set:Nn \l_delimiter_size_tl {}
  }{}
  \pgfmathparse{\@nargin==0 ? 1:0}
  \ifthenelse{\pgfmathresult>0}{
    \tl_set:Nn \l_norm_type_tl {}
  }{}
  \ifthenelse{\equal{\l_delimiter_size_tl}{auto}}{%
    \@@normgen[auto]{#2}}{%
    \ifthenelse{\equal{\l_delimiter_size_tl}{}}{
      \@@normgen[]{#2}}{\@@normgen[\l_delimiter_size_tl]{#2}}}%
  \c_math_subscript_token{\l_norm_type_tl}
  \group_end:
}
\seq_new:N \l_var_seq
\tl_new:N \l_norm_type_tl
\tl_new:N \l_delimiter_size_tl
\makeatother
\ExplSyntaxOff

\newcommand{\vertiii}[1]{{\left\vert\kern-0.25ex\left\vert\kern-0.25ex\left\vert #1
        \right\vert\kern-0.25ex\right\vert\kern-0.25ex\right\vert}}

\newcommand{\seq}[1]{\{#1\}}

\newlength{\dhatheight}

\makeatletter
\define@key{projKeys}{set}{\def\set{#1}}
\define@key{proxKeys}{t}{\def\t{#1}}
\define@key{proxKeys}{f}{\def\f{#1}}
\makeatother

\DeclareMathOperator{\proximal}{prox}

\newcommand{\Prox}[2][]{%
  \setkeys{proxKeys}{t=t,#1}%
  \setkeys{proxKeys}{f=f,#1}%
  \mathchoice{\underset{\t\f}{\proximal}}%
  {\proximal_{\t\f}}{\proximal_{\t\f}}{\proximal_{\t\f}}%
  \ifthenelse{\equal{#2}{}}{}{\left(#2\right)}%
}
\DeclareMathOperator{\distance}{dist}
\newcommand{\dist}[2][]{%
  \distance_{#1}\ifthenelse{\equal{#2}{}}{}{(#2)}
}
\newcommand{\gauge}[2][]{%
  \gamma_{#1}\ifthenelse{\equal{#2}{}}{}{(#2)}
}

\newcommand{\transp}{{\scriptscriptstyle\mathsf{T}}}
\newcommand{\T}{^\transp}

\newcommand{\dd}[1][]{#1\mathrm{d}}

\newcommand{\fun}[2][1]{%
  #2(%
  \foreach \index in {1, ..., #1} {%
    \ifthenelse{\equal{\index}{#1}}{%
      \cdot%
    }{%
      \cdot,%
    }%
  })}
\newcommand{\definedas}[1][tri]{%
  \ifthenelse{\equal{#1}{tri}}{\triangleq}{\coloneqq}}

\renewcommand{\implies}{\Rightarrow}

\DeclareMathOperator*{\exptx}{exp}
\renewcommand{\exp}[2][exponent]{\ifthenelse{\equal{#1}{exponent}}{e^{#2}}{\exptx\left(#2\right)}}
\DeclareMathOperator*{\expsf}{\mathsf{exp}}
\newcommand{\expm}[2][exponent]{\ifthenelse{\equal{#1}{exponent}}{\mathsf{e}^{#2}}{\expsf\left(#2\right)}}

\makeatletter
\DeclareFontFamily{U}{tipa}{}
\DeclareFontShape{U}{tipa}{m}{n}{<->tipa10}{}
\newcommand{\arc@char}{{\usefont{U}{tipa}{m}{n}\symbol{62}}}%
\renewcommand{\arc}[1]{\mathpalette\arc@arc{#1}}
\newcommand{\arc@arc}[2]{%
  \sbox0{$\m@th#1#2$}%
  \vbox{
    \hbox{\resizebox{\wd0}{\height}{\arc@char}}
    \nointerlineskip
    \box0
  }%
}
\makeatother



\newcommand{\dash}{\Hyphdash}

\makeatletter
\newcommand{\pushright}[1]{%
  \ifmeasuring@#1\else\omit\hfill$\displaystyle#1$\fi\ignorespaces}
\newcommand{\pushleft}[1]{%
  \ifmeasuring@#1\else\omit$\displaystyle#1$\hfill\fi\ignorespaces}
\makeatother

\newcommand{\union}{\cup}

\newcommand{\set}[1]{\mathcal #1}

\newcommand{\eucledian}[1][]{\ifthenelse{\equal{#1}{}}{\mathbb E}{\mathbb #1}}
\DeclareMathOperator{\@convhull}{conv}
\newcommand{\reals}[1][]{%
  \ifthenelse{\equal{#1}{extended}}{\overline}{}{\mathbb R}%
}
\newcommand{\real}{\mathbb R}

\newcommand{\quaternion}{\mathbb Q}

\newcommand{\algvar}[1]{\text{\IfSubStr{#1}{_}{%
    \StrSubstitute{#1}{_}{\textunderscore}}{#1}}}

\makeatletter
\definecolor{listinggray}{gray}{0.9}
\definecolor{lbcolor}{rgb}{0.9,0.9,0.9}
\colorlet{Darkgreen}{green!60!black}
\lstdefinelanguage{Julia}%
{morekeywords={abstract,break,case,catch,const,continue,do,else,elseif,%
    end,export,false,for,function,immutable,import,importall,if,in,%
    macro,module,otherwise,quote,return,switch,true,try,type,typealias,%
    using,while},%
  sensitive=true,%
  alsoother={\$},
  morecomment=[l]\#,%
  morecomment=[n]{\#=}{=\#},%
  morestring=[s]{"}{"},%
  morestring=[m]{'}{'},%
}[keywords,comments,strings]
\lstdefinestyle{listing@C++}{
  language=[GNU]C++,
}
\lstdefinestyle{listing@Julia}{%
  language = Julia,
}
\lstdefinestyle{listing@Python}{%
  language = Python,
}
\newcommand\langname@bash{}
\def\langname@bash{bash}
\newcommand\prompt@bash{\texttt{\$}\ } 
\newcommand\addedToEveryPar@bash{}
\lst@AddToHook{EveryPar}{\addedToEveryPar@bash}
\lst@AddToHook{PreInit}{%
  \ifx\lst@language\langname@bash%
  \let\addedToEveryPar@bash\prompt@bash%
  \fi
}
\lstdefinestyle{listingterminal}{%
  language=bash,
  backgroundcolor=\color{black!90!blue},
  basicstyle=\color{white},
  commentstyle=\color{yellow!50},
  keywordstyle=\bfseries\color{orange},
  breaklines,
  escapechar=@
}

\lstnewenvironment{terminalcode}[1][]{\lstset{style=listingterminal,#1}}{}
\makeatother

\newcommand{\subdiff}[1][none]{%
  \ifthenelse{\equal{#1}{none}}{%
    \partial%
  }{%
    \partial_{#1}%
  }%
}

\makeatletter
\DeclareMathOperator{\@convhull}{conv}
\newcommand{\convhull}[1][]{%
  \ifthenelse{\equal{#1}{closed}}{\overline}{}\@convhull%
}
\DeclareMathOperator{\@closed@convex@envelope}{\overline{co}}
\newcommand{\cce}{\@closed@convex@envelope}
\makeatother

\newcommand{\support}[2][delta]{\ifthenelse{\equal{#1}{delta}}{%
    \delta^*_{#2}%
  }{\sigma_{#2}}}
\newcommand{\conj}{^{*}}

\newcommand{\Matrix}[2][]{%
  \ifthenelse{\equal{#1}{}}{}{\setlength\arraycolsep{#1}}%
  \begin{bmatrix}#2\end{bmatrix}}
\renewcommand{\Array}[2][]{%
  \ifthenelse{\equal{#1}{}}{}{\setlength\arraycolsep{#1}}%
  \begin{matrix}#2\end{matrix}}

\makeatletter
\newcommand{\toset}{%
  \def\arr@offset{0.15em}
  \def\arr@len{0.7em}
  \def\arr@height{0.3em}
  \tikz[minimum height=0ex,outer sep=0,inner sep=0]
  \path[-{Latex[length=0.8mm]}]
  node (a) at (0,0) {}
  node (b) at (0,\arr@height) {}
  (a) edge ++(\arr@len,0)
  (b) edge ++(\arr@len,0)
  (a) edge[draw=none] ++(0,-\arr@offset);%
}
\makeatother

\newcommand{\ev}[3][c]{%
  \ifthenelse{\equal{#1}{c}}{%
    \ifthenelse{\equal{#2}{}}{\forall[0,#3]}{\forall[#2,#3]}%
  }{%
    \ifthenelse{\equal{#1}{o}}{%
      \ifthenelse{\equal{#2}{}}{\forall(0,#3)}{\forall(#2,#3)}%
    }{%
      \ifthenelse{\equal{#1}{oc}}{%
        \ifthenelse{\equal{#2}{}}{\forall(0,#3]}{\forall(#2,#3]}%
      }{%
        \ifthenelse{\equal{#2}{}}{\forall[0,#3)}{\forall[#2,#3)}%
      }%
    }%
  }%
}

\DeclareMathOperator*{\argmax}{argmax}

\DeclareMathOperator{\val}{val}

\newcommand*\softmax{\lse}

\newcommand{\GetLabel}[1]{\expandafter\csname #1Label\endcsname}
\makeatletter
\define@key{printRefKeys}{otherLabel}{\def\otherLabel{#1}}
\define@key{printRefKeys}{concatenate}{\def\concatenate{#1}}
\makeatother
\newcommand{\PrintRefs}[3][]{%
  \setkeys{printRefKeys}{otherLabel=,#1}%
  \setkeys{printRefKeys}{concatenate=false,#1}%
  \xdef\MyModifiedLabel{#2}%
  \ifthenelse{\equal{\otherLabel}{}}{}{\xdef\MyModifiedLabel{\otherLabel}}%
  \xdef\MyRefCount{0}%
  \foreach \i in {#3} {%
    \tikzmath{\MyRefCount=int(\MyRefCount+1);}%
    \xdef\MyRefCount{\MyRefCount}%
  }%
  \edef\OrigRefCount{\MyRefCount}%
  \ifthenelse{\equal{\MyRefCount}{1}}{%
    #2~\ref{\GetLabel{\MyModifiedLabel}:#3}%
  }{%
    \xdef\MyCounter{0}%
    \ifthenelse{\equal{#2}{Corollary}}{%
      Corollaries%
    }{%
      #2s%
    }~%
    \foreach \AlgRef in {#3} {%
      \tikzmath{\MyCounter=int(\MyCounter+1);}%
      \xdef\MyCounter{\MyCounter}%
      \ifthenelse{\equal{\concatenate}{true}}{%
        \pgfmathparse{\MyCounter==1 ? 1 : 0}%
        \ifthenelse{\pgfmathresult>0}{%
          \ref{\GetLabel{\MyModifiedLabel}:\AlgRef}-%
        }{%
          \pgfmathparse{\MyCounter==\MyRefCount ? 1 : 0}%
          \ifthenelse{\pgfmathresult>0}{%
            \ref{\GetLabel{\MyModifiedLabel}:\AlgRef}%
          }{}%
        }%
      }{%
        \pgfmathparse{\MyCounter<\MyRefCount-1 ? 1 : 0}%
        \ifthenelse{\pgfmathresult>0}{%
          \ref{\GetLabel{\MyModifiedLabel}:\AlgRef}, %
        }{%
          \pgfmathparse{\MyCounter<\MyRefCount ? 1 : 0}%
          \ifthenelse{\pgfmathresult>0}{%
            \ref{\GetLabel{\MyModifiedLabel}:\AlgRef}%
            \pgfmathparse{\OrigRefCount==2 ? 1:0}%
            \ifthenelse{\pgfmathresult>0}{ and }{, and }%
          }{%
            \ref{\GetLabel{\MyModifiedLabel}:\AlgRef}%
          }%
        }%
      }%
    }%
  }%
}
\newcommand{\sref}[1]{\PrintRefs{Section}{#1}}
\newcommand{\ssref}[2][]{\PrintRefs[otherLabel=Subsection,#1]{Section}{#2}}
\newcommand{\sssref}[2][]{\PrintRefs[otherLabel=Subsubsection,#1]{Section}{#2}}
\newcommand{\cref}[1]{\PrintRefs{Chapter}{#1}}

\newcommand{\figref}[1]{\PrintRefs{Figure}{#1}}
\newcommand{\tabref}[2][]{\PrintRefs[#1]{Table}{#2}}

\makeatletter
\define@key{algKeys}{start}{\def\startline{#1}}
\define@key{algKeys}{end}{\def\endline{#1}}
\define@key{algKeys}{show}{\def\showalg{#1}}
\makeatother
\renewcommand{\algref}[2][]{%
  \setkeys{algKeys}{start=,#1}%
  \setkeys{algKeys}{end=,#1}%
  \setkeys{algKeys}{show=true,#1}%
  \ifthenelse{\equal{\startline}{}}{}{line%
    \ifthenelse{\equal{\endline}{}}{}{s}}%
  \ifthenelse{\equal{\startline}{}}{}{~\ref{alg:#2:line:\startline}%
    \ifthenelse{\equal{\endline}{}}{}{-\ref{alg:#2:line:\endline}} %
    of }%
  \ifthenelse{\equal{\showalg}{true}}{\PrintRefs{Algorithm}{#2}}{}%
}
\newcommand{\pref}[2][]{\PrintRefs[#1]{Problem}{#2}}

\makeatletter
\newcommand*{\addFileDependency}[1]{
  \typeout{(#1)}
  \@addtofilelist{#1}
  \IfFileExists{#1}{}{\typeout{No file #1.}}
}
\makeatother

\DeclareSIUnit{\radian}{rad}
\DeclareSIUnit\century{century}
\DeclareSIUnit\year{yr}
\DeclareSIUnit{\ton}{t}

\makeatletter

\newcommand{\add@list@item}[2]{
  \ifthenelse{\equal{#1}{}}{\xdef#1{{#2}}}{\xdef#1{#1,{#2}}}
}

\def\providecounter#1{%
  \@ifundefined{c@#1}%
  {\newcounter{#1}}{\setcounter{#1}{0}}}

\newcommand{\makelist}[2]{
  \xdef\@current@list@name{#1}
  \xdef\@current@list@counter{#1@listcounter}
  \providecounter{\@current@list@counter}
  \edef\@tmp@list{{#2}}
  \expandafter\forcsvlist\expandafter\list@saveitem\@tmp@list
}

\newcommand{\list@saveitem}[1]{%
  \stepcounter{\@current@list@counter}%
  \expandafter\def\csname\@current@list@name%
  \arabic{\@current@list@counter}\endcsname{#1}
}

\newcommand{\list@nth}[2]{\csname #1#2\endcsname}

\makeatother

\newcommand{\mrm}[2][m]{\ifthenelse{\equal{#1}{m}}{\mathrm{#2}}
  {\textnormal{#2}}}

\newcommand{\Behcet}{Beh\c{c}et}
\newcommand{\Acikmese}{A\c{c}{\i}kme\c{s}e}

\makeatletter
\@ifclassloaded{beamer}{}{}
\makeatother

\newcommand{\makeflag}[2]{%
  \expandafter\newif\csname ifmake#1\endcsname
  \csname make#1#2\endcsname
}

\NewEnviron{makeif}[2][]{%
  \csname ifmake#2\endcsname
  \BODY
  \else
  #1
  \fi
}

\ifthenelse{\isundefined{\noindentadjust}}{
\setlength\parindent{0pt} 
\setlength{\parskip}{4pt} 
}{}
\allowdisplaybreaks       

\makeatletter
\ifthenelse{\isundefined{\nohref}}{
\@ifclassloaded{beamer}{}{
\ifthenelse{\isundefined{\preprintcsm}}{
  \usepackage[hypertexnames=false]{hyperref}}{}}{}}{}
\makeatother

\ifthenelse{\isundefined{\nohref}}{%
  \hypersetup{colorlinks,linkcolor={red},citecolor={green},urlcolor={blue}}
  }{}

\newtheoremstyle{faded}
{\topsep}%
{\topsep}%
{\color{black!75!white}}%
{}%
{\color{black}\bfseries}
{.}
{.5em}%
{\thmname{#1}~\thmnumber{#2}\thmnote{ (#3)}}%

\ifthenelse{\isundefined{\theoremlook}}{
\def\theoremlook{theorem}}{}

\ifthenelse{\isundefined{\definitionlook}}{
\def\definitionlook{definition}}{}

\theoremstyle{\theoremlook}

\theoremstyle{\definitionlook}

\makeatletter
\def\thmenvs{theorem,lemma,corollary,proposition,remark,property,%
  condition,method,example,problem,definition,assumption}
\foreach \@te in \thmenvs {
  \bgroup
  \globaldefs=1
  \expandafter\let\csname o\@te\expandafter\endcsname\csname\@te\endcsname
  \expandafter\let\csname eo\@te\expandafter\endcsname\csname end\@te\endcsname
  \renewenvironment{\@te}[1][]{%
    \csname o\@currenvir\endcsname%
    \ifthenelse{\equal{##1}{}}{}{\label{\@currenvir:##1}}
  }{%
    \expandafter\csname eo\@currenvir\endcsname%
  }
  \egroup
}
\makeatother

\newcommand{\defvar}[3][show]{\nomenclature{#2}{#3}%
  \ifthenelse{\equal{#1}{show}}{#2}{}}

\input{latex_common/optimization.tex}
\input{latex_common/review_commands.tex}
\hypersetup{colorlinks,linkcolor={red},citecolor={green},urlcolor={blue}}

\definecolor{DarkBlue}{HTML}{26415d}
\definecolor{Blue}{HTML}{356397}
\definecolor{Green}{HTML}{5da9a1}
\definecolor{Yellow}{HTML}{f1d46a}
\definecolor{Red}{HTML}{db6245}
\definecolor{Black}{HTML}{000000}
\definecolor{Grey}{HTML}{acacac}

\renewcommand{\frame}[1]{\mathcal F_{#1}}

\newcommand{\tf}{t_f}
\newcommand{\nrcs}{{n_{\mrm{rcs}}}}
\newcommand{\Frcs}{F_{\mrm{rcs}}}
\newcommand{\Dt}[1][]{\Delta t_{#1}}
\newcommand{\Dtr}[1][]{\Delta t'_{#1}}
\newcommand{\norb}{n_o}
\newcommand{\ex}[1][]{\hat x_{#1}}
\newcommand{\ey}[1][]{\hat y_{#1}}
\newcommand{\ez}[1][]{\hat z_{#1}}
\newcommand{\tc}{t_c}
\newcommand{\Nc}{N_c}
\newcommand{\Ifr}{\mathcal I_{\mrm{fr}}}
\newcommand{\fdir}[1][]{\hat f_{#1}}
\renewcommand{\th}[1]{#1-th}
\newcommand{\rplume}{r_{\mrm{plume}}}
\newcommand{\zplume}{\zeta_{\mrm{plume}}}
\newcommand{\Mplume}{M_{\mrm{plume}}}
\newcommand{\rappch}{r_{\mrm{appch}}}
\newcommand{\angappch}{\theta_{\mrm{appch}}}
\newcommand{\qL}{q_{\ell}}
\newcommand{\qDP}{q_{\mrm{dp}}}
\newcommand{\rDP}{p_{\mrm{dp}}}
\newcommand{\Relax}[1]{\Delta #1}
\newcommand{\Dxf}{\Relax{x_f}}
\newcommand{\tol}[1]{\varepsilon_{#1}}
\newcommand{\Jfuel}{J_{\mrm{fuel}}}
\newcommand{\Jeq}{J_{\mrm{eq}}}
\newcommand{\weq}{w_{\mrm{eq}}}
\newcommand{\minusone}[1]{#1'}
\renewcommand{\ng}{n_g}
\newcommand{\nz}{n_z}
\newcommand{\nf}{n_f}
\newcommand{\predicate}[1][]{g_{#1}}
\newcommand{\Lif}{L}
\newcommand{\Rif}{R}
\newcommand{\fL}{f_{\Lif}}
\newcommand{\fR}{f_{\Rif}}
\renewcommand{\sharp}[1][]{\kappa_{#1}}
\newcommand{\sigmoid}[1][]{\sigma_{#1}}
\newcommand{\smooth}[1]{\tilde{#1}_{\sharp}}

\renewcommand{\softmax}{\mathtt{LSE}_{\sharp}}
\newcommand{\predmatch}{{\predicate}_c}
\newcommand{\textcode}[1]{\texttt{#1}\xspace}
\newcommand{\texttrue}{\textcode{true}}
\newcommand{\textfalse}{\textcode{false}}
\newcommand{\textand}{\textcode{and}}
\newcommand{\textor}{\textcode{or}}
\newcommand{\textif}{\textcode{if}}
\newcommand{\textelse}{\textcode{else}}
\newcommand{\textifelse}{\textcode{if-else}}
\newcommand{\Dtbuff}{\Dt[\mrm{db}]}
\newcommand{\Gbuff}{G_{\mrm{db},\sharp}}
\newcommand{\smoothRappch}{\smooth{\Rif}^{\mrm{appch}}}
\newcommand{\smoothRplume}{\smooth{\Rif}^{\mrm{plume}}}
\newcommand{\smoothRmib}{\smooth{\Rif}^{\mrm{mib}}}
\newcommand{\gic}{g_{\mrm{ic}}}
\newcommand{\gtc}{g_{\mrm{tc}}}
\newcommand{\dimx}{n_x}
\newcommand{\dimu}{n_u}
\newcommand{\dimp}{n_p}
\newcommand{\dimss}{n_s}
\newcommand{\dimgic}{n_{\mrm{ic}}}
\newcommand{\dimgtc}{n_{\mrm{tc}}}
\newcommand{\Nhom}{N_h}

\newcommand{\solution}[1][]{X^{*}_{\ifthenelse{\equal{#1}{}}{}{(#1)}}}
\newcommand{\interim}[1][]{X_{\ifthenelse{\equal{#1}{}}{}{(#1)}}}
\newcommand{\Li}{L}
\newcommand{\li}{\ell}
\newcommand{\sigeps}{\varepsilon}
\newcommand{\sigdelta}{\delta}
\newcommand{\sigdeltamin}{\sigdelta_1}
\newcommand{\sigdeltamax}{\sigdelta_0}
\newcommand{\siginterp}{\alpha}
\newcommand{\sigcoeff}[1][]{\rho_{#1}}
\newcommand{\cost}[1][]{J_{#1}}
\newcommand{\hombl}{\beta_{\mrm{worse}}}
\newcommand{\hombu}{\beta_{\mrm{trig}}}

\newcommand{\frL}{\mathcal L} 
\newcommand{\frB}{\mathcal B} 

\newcommand{\alglocation}[1]{%
  \tikz[baseline=-0.6ex]{\node[
    circle,
    draw=black,
    minimum size=0.6em,
    inner sep=0.2mm,
    scale=0.8
    ] at (0,0) {#1};}}

\def\iterstartloc{1}
\def\solveloc{2}
\def\testloc{3}

\newcommand{\resultplot}[1]{code/ptr_rendezvous_3d_#1.pdf}

\renewcommand*{\arraystretch}{0.8}

\graphicspath{{./figures/}}

\title{Fast Homotopy for Spacecraft Rendezvous Trajectory Optimization with
  Discrete Logic}

\author{%
  Danylo Malyuta\footnote{Ph.D. Candidate, Autonomous Controls Laboratory,
    Aeronautics \& Astronautics, AIAA Student Member. Contact:
    \texttt{danylo@uw.edu}.}
  and \Behcet{} \Acikmese{}\footnote{Professor, Autonomous Controls Laboratory,
    Aeronautics \& Astronautics, AIAA Associate Fellow. Contact:
    \texttt{behcet@uw.edu}.}}

\affil{University of Washington, Seattle, Washington 98195}

\begin{document}

\maketitle

\begin{abstract}

  This paper presents a computationally efficient optimization algorithm for
  solving nonconvex optimal control problems that involve discrete logic
  constraints. Traditional solution methods for these constraints require
  binary variables and mixed\dash integer programming, which is prohibitively
  slow and computationally expensive. This paper targets a fast solution that
  is capable of real\dash time implementation onboard spacecraft. To do so, a
  novel algorithm is developed that blends sequential convex programming and
  numerical continuation into a single iterative solution process. Inside the
  algorithm, discrete logic constraints are approximated by smooth functions,
  and a homotopy parameter governs the accuracy of this approximation. As the
  algorithm converges, the homotopy parameter is updated such that the smooth
  approximations enforce the exact discrete logic. The effectiveness of this
  approach is numerically demonstrated for a realistic rendezvous scenario
  inspired by the Apollo Transposition and Docking maneuver. In under 15
  seconds of cumulative solver time, the algorithm is able to reliably find
  difficult fuel\dash optimal trajectories that obey the following discrete
  logic constraints: thruster minimum impulse\dash bit, range\dash triggered
  approach cone, and range\dash triggered plume impingement. The optimized
  trajectory uses significantly less fuel than reported NASA design targets.

\end{abstract}

\makenomenclature
\printnomenclature

\section{Introduction}
\label{section:intro}

\lettrine{S}{pace} programs have historically been deemed mature once they
establish the ability to perform rendezvous and docking
operations~\cite{Woffinden2007}. Some of the earliest programs of the United
States and the Soviet Union (e.g., Gemini and Soyuz) had as their explicit goal
to demonstrate the capability of performing rendezvous, proximity operations,
and docking maneuvers. The ultimate objective to land humans on the moon drove
the need for these capabilities. Beyond the lunar missions of the 1960s,
rendezvous and docking continued to be a core technology required to construct
and service space stations that were built in low Earth
orbit~\cite{Goodman2006}. The Shuttle program was comprised of dozens of
missions for which rendezvous (and more generally, proximity operations) was an
explicit mission objective. The core technology used to achieve rendezvous and
docking has remained largely unchanged in the decades since the earliest
maneuvers were successful. While this heritage technology is far from obsolete,
it has been stated that it may be unable to meet the requirements of future
missions~\cite{Woffinden2007}. A driving force that will require new methods is
the need for a system that can perform \textit{fully autonomous} rendezvous in
several domains (e.g., low Earth orbit, low lunar orbit,
etc.)~\cite{DSouza2007}. Several vehicles capable of autonomous docking are
either already operational or in development, ranging from large vehicles such
as the SpaceX Crew Dragon, Soyuz, and Orion
\cite{Woffinden2007,DSouza2007,SpaceXAutoDock}, to smaller robotic vehicles for
clearing orbital debris \cite{Nishida2011,Kaplan2009,SpaceJunkRendezvous}.

The objective of this paper is to present a framework for designing autonomous
docking trajectories that accurately reflect the capabilities and constraints
that have been historically prevalent for proximity operation missions. We view
the problem as a trajectory generation problem, and compute what would be
implemented as a guidance solution. In particular, we show how to model
challenging \textit{discrete logic} constraints within a continuous
optimization framework. The resulting algorithm is numerically demonstrated to
be sufficiently fast for ground\dash based use, and has the potential to be
real\dash time capable if implemented in a compiled programming language. A
link to the numerical implementation of the algorithm is provided in
\cite{OurOpenSourceCode}.

The open-loop generation of spacecraft docking trajectories using
optimization-based methods is a relatively new field spawned by the shift
towards autonomy \cite{CSM2021}. Open\dash loop trajectory generation computes a
complete start\dash to\dash finish trajectory, and leaves robust tracking to
closed\dash loop feedback control. 
In~\cite{Miele2007,Miele2007a} the authors discuss both time- and fuel\dash
optimal solutions with a focus on problem formulations that are conducive to
on-board implementation. Their study offers an insightful view on the structure
of optimality at the cost of a simplified problem statement and omission of
state constraints. In \cite{Pascucci2017}, lossless convexification is used to
generate fuel\dash optimal docking trajectories which account for nonconvex
thrust and plume impingement constraints, albeit the thrust is not allowed to
turn off. In \cite{Harris2014b}, lossless convexification allows to generate
bang-bang controls for minimum-time spacecraft rendezvous using differential
drag, however without state constraints or spacecraft attitude dynamics. A
similar relaxation is also presented in \cite{Lu2013}, where a sequential convex
programming (SCP) algorithm is developed for near\dash field autonomous
rendezvous in an arbitrary Keplerian orbit. Range\dash triggered approach cone
and plume impingement constraints are imposed, however their activation is
pre-determined through user specification rather than automatically by the
algorithm. A similar solution method is considered in \cite{Liu2013}, where a
rendezvous problem is solved with aerodynamic drag, gravity harmonics, and a
nonconvex keep\dash out ellipsoid for collision avoidance. The latter
constraint applies during the initial maneuvering phase, while for the final
approach the keep\dash out zone is replaced by a convex approach cone.

In \cite{Breger2008}, an optimization framework is used to impose safety-based
constraints in the case of anomalous behavior (including thruster failure) by
introducing a suboptimal convex program to design safe trajectories which
approximate a nonconvex mixed-integer problem using a new set of ``safe''
inputs. Along the same lines of mixed-integer programming,~\cite{Richards2002}
solves a fuel\dash optimal problem subject to thrust plume and collision
avoidance constraints. The authors introduce several heuristic techniques in
order to fit the problem within the scope of mixed-integer linear programming,
but still observe rather long solve times (over 40 minutes in some cases). More
recently,~\cite{Sun2019} studied a multi-phase docking problem with several
state constraints. The authors use binary variables to impose different
constraints during each phase, and propose an iterative solution method with
closed-form update rules. Beyond the use of mixed-integer
methods,~\cite{Phillips2003} proposes a randomized optimization method similar
to the $A^*$ method, while~\cite{Hartley2013} proposes a convex one\dash norm
regularized model predictive control solution.

Notably, the aforementioned references do not consider the spacecraft attitude
during trajectory generation and do not explicitly account for what is referred
to as the minimum impulse\dash bit (MIB) of the reaction control thrusters that
are used to realize the trajectories. The latter constraint refers to the fact
that impulsive chemical thrusters cannot fire for an arbitrarily short
duration, since there is some minimum pulse width that is inherent to the
hardware. Hartley et al. \cite{Hartley2013} acknowledge this issue, but instead
of explicitly enforcing the constraint, the authors use a one\dash norm penalty
term to discourage violation of the constraint (i.e., a soft constraint). Our
view is that both attitude and the MIB constraint are critical for close
proximity operations such as the terminal phase of rendezvous and docking,
where two spacecraft are maneuvering close to each other. We thus target an
algorithm that can efficiently incorporate both effects.

\subsection{Contributions}
\label{subsection:intro:contributions}

This paper's contribution is a numerical optimization algorithm to solve
optimal control problems (OCPs) that involve a general class of discrete logic
constraints. The algorithm is based on a novel arrangement of two core
methodologies: sequential convex programming and numerical continuation. SCP is
a trust region method for solving general nonconvex optimal control problems
\cite{CSM2021}. However, it is incapable of handling discrete constraints in
their pure (integer) form. By using a homotopy map based on the multinomial
logit function, we embed smooth approximations of discrete constraints into the
SCP framework, a process also known as continuous embedding
\cite{Bengea2005}. The homotopy map is then updated via a numerical continuation
scheme, which transforms an initial coarse approximation into an arbitrarily
precise representation of the discrete logic. Herein lies our key innovation:
we run SCP and numerical continuation \textit{in parallel}, rather than in the
traditional sequenced approach where one homotopy update is followed by a full
SCP solve. For this reason, we call the method \textit{embedded} numerical
continuation. The resulting algorithm is shown to converge quickly and reliably
for a representative terminal rendezvous problem inspired by the Apollo
Transposition and Docking maneuver. The problem involves the following major
constraints: full six degree of freedom (DOF) dynamics, thruster minimum
impulse\dash bit, range\dash triggered approach cone, and range\dash triggered
plume impingement. The latter constraints are similar to those considered in
\cite{Lu2013,Liu2013}, with the advantage that discrete logic allows the
approach cone and plume impingement constraints to be switched on/off
\textit{automatically} by the algorithm, without user input.

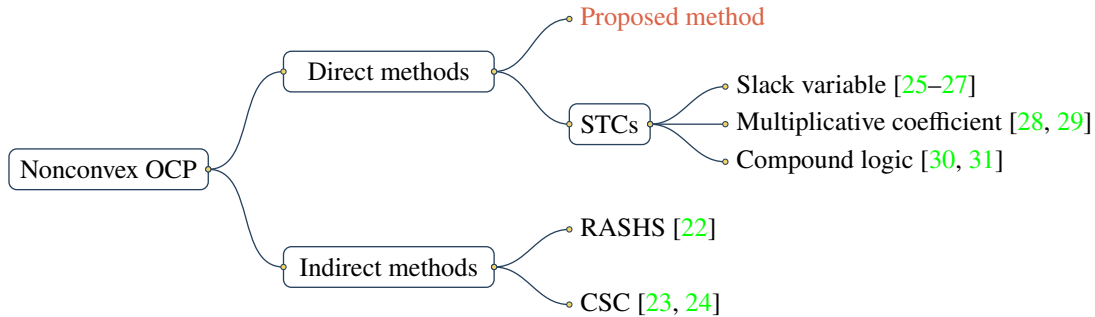
\begin{figure}
  \centering




\makeatletter

\def\hsep{10mm}
\def\method@vsep{13mm}
\def\algo@vsep{7mm}
\def\subalgo@vsep{5mm}

\tikzset{
  box/.style={
    line width=0.5pt,
    draw=DarkBlue,
    rounded corners=3pt,
    inner sep=1.5mm,
    anchor=west
  },
  box noborder/.style={
    box,
    draw=none
  },
  connector line/.style={
    line width=0.5pt,
    draw=DarkBlue
  },
  connector endpoint/.style={
    circle,
    fill=Yellow,
    draw=DarkBlue,
    minimum size=2pt
  }
}

\newcommand{\connect}[2]{
  \draw[connector line] (#1) to[out=0, in=180]
  node[pos=0](connector start){}
  node[pos=1](connector end){} (#2);
  \node[connector endpoint] at (connector start) {};
  \node[connector endpoint] at (connector end) {};
}

\begin{tikzpicture}[line cap=round, line join=round]

  \coordinate (origin) at (0, 0);


  \node[box] (ocp) at (origin) {Nonconvex OCP};


  \node[box,shift={(\hsep,\method@vsep)},
  minimum width=2.8cm] (direct) at (ocp.east)
  {Direct methods};

  \node[box,shift={(\hsep,-\method@vsep)},
  minimum width=2.8cm] (indirect) at (ocp.east)
  {Indirect methods};

  \connect{ocp.east}{direct.west};
  \connect{ocp.east}{indirect.west};


  \node[box noborder,shift={(\hsep,\subalgo@vsep)}] (rashs) at (indirect.east)
  {RASHS \cite{saranathan2018relaxed}};
  \node[box noborder,shift={(\hsep,-\subalgo@vsep)}] (csc) at (indirect.east)
  {CSC \cite{taheri2020novel,arya2021composite}};

  \connect{indirect.east}{rashs.west};
  \connect{indirect.east}{csc.west};


  \node[box,shift={(\hsep,-\algo@vsep)}] (stcs) at (direct.east)
  {STCs};

  \node[box noborder,shift={(\hsep,\subalgo@vsep)}] (miki) at (stcs.east)
  {Slack variable \cite{Szmuk2020,MalyutaScitechDocking,SzmukThesis}};
  \node[box noborder,shift={(\hsep,0)}] (taylor) at (stcs.east)
  {Multiplicative coefficient \cite{Reynolds2020,ReynoldsThesis}};
  \node[box noborder,shift={(\hsep,-\subalgo@vsep)}] (compound) at (stcs.east)
  {Compound logic \cite{Szmuk2019a,Szmuk2019b}};

  \connect{direct.east}{stcs.west}
  \connect{stcs.east}{miki.west}
  \connect{stcs.east}{taylor.west}
  \connect{stcs.east}{compound.west}


  \node[box noborder,shift={(\hsep,\algo@vsep)},
  text=Red] (us) at (direct.east)
  {Proposed method};

  \connect{direct.east}{us.west}

\end{tikzpicture}

\makeatother


  \caption{Illustration of the proposed algorithm's relationship to existing
    literature on handling discrete logic in a continuous-variable optimization
    framework.}
  \label{fig:algorithm_tree}
\end{figure}

This paper represents a significant upgrade in terms of both runtime
performance and convergence reliability over the same authors' previous
publication on SCP\dash based rendezvous
\cite{MalyutaScitechDocking}. \figref{algorithm_tree} illustrates how the proposed
algorithm relates to existing literature. Closest to our approach are the
recently published relaxed autonomous switched hybrid system (RASHS) and
composite smooth control (CSC) algorithms
\cite{saranathan2018relaxed,taheri2020novel,arya2021composite}. Both RASHS and
CSC belong to the indirect family of methods in the sense that they solve for
the optimality conditions obtained from Pontryagin's maximum principle
\cite{Betts1998,PontryaginBook,BerkovitzBook}. Furthermore, both RASHS and CSC
handle discrete logic that is representable by a sequence of Boolean \textand
gates. Our method is distinct from RASHS and CSC in two ways. First, it is a
direct method in the sense that it uses numerical optimization to solve a
discretized version of the optimal control problem. This generally makes it
easier to handle constraints, which are nontrivial to include in an indirect
approach. Second, the proposed method models discrete logic that is
representable by a sequence of Boolean \textor gates. As a result, our work
combined with RASHS and CSC can extend homotopy to general Boolean logic using
any combination of logic gates. A more detailed comparison of the methods is
given in \sssref{embedding:description:existing}.

Our algorithm is also closely related to the recently introduced family of
state triggered constraints (STCs) for SCP algorithms
\cite{Szmuk2020,Reynolds2020}. Unlike our method, STCs directly use
linearization instead of homotopy in order to enforce an equivalent
continuous-variable formulation of discrete logic constraints. Several versions
of STCs have been introduced, and we cover these in more detail in
\ssref{embedding:motivation}. Past work on STCs, however, discovered that they can
exhibit unfavorable ``locking'' behavior for thruster minimum impulse-bit
constraints that are relevant for spacecraft rendezvous
\cite{MalyutaScitechDocking}. This phenomenon prevents the algorithm from
converging, and we describe it in detail in \ssref{embedding:motivation}. The
algorithm presented in this article handles discrete logic constraints like
STCs, and does not exhibit locking.


\subsection{Structure}
\label{subsection:intro:outline}

The rest of this paper is structured as follows. In \sref{formulation} we formulate
the rendezvous problem that is to be solved, but which is not efficiently
solvable in its raw form. \sref{embedding} then describes the homotopy map which
can model a generic class of discrete logic in a smooth way. Using this
smoothing, \sref{scp} describes our key contribution: an algorithm that can solve
nonconvex optimal control problems with discrete logic. The effectiveness of
the approach is numerically demonstrated in \sref{results} for a realistic scenario
based on the historical Apollo Transposition and Docking maneuver.

The paper uses mostly standard mathematical notation. However, it is worth
emphasizing the following specific notational elements. Inline vector
concatenation is written as $\brak[big]{a;b;c}\in\reals^3$ where
$a,b,c\in\reals$. To avoid confusion, faded text is used to denote expressions
that belong to a summation, such as $\sum_{i=1}^n \fade{x_i}+z$ where $z$ is
outside the sum. The symbol $\union$ stands for set union, while the symbols
$\land$ and $\lor$ represent Boolean \textand and \textor operators. Quaternion
multiplication is denoted by $\otimes$.









\section{Rendezvous Problem Formulation}
\label{section:formulation}

In this section we formulate a trajectory generation problem where the
objective is to guide a chaser spacecraft to dock with a passive target
spacecraft in a predetermined orbit. We assume that the maneuver happens in low
Earth orbit (LEO) and that the target's orbit is circular. The chaser's
dynamics are defined in \ssref{formulation:chaser_dynamics}, the actuator model is
described in \ssref{formulation:impulsive_thrust}, and the rendezvous constraints
are given in \ssref{formulation:plume_impingement, formulation:approach_cone,
  formulation:bcs}. \ssref{formulation:summary} gives a complete formulation of the
free\dash final time nonconvex optimal control problem which, if solved,
generates a fuel\dash optimal rendezvous trajectory. Most notably, because the
constraints in \ssref{formulation:impulsive_thrust, formulation:plume_impingement,
  formulation:approach_cone} involve discrete logic, the problem is not readily
solvable by traditional continuous optimization methods.

\subsection{Chaser Spacecraft Dynamics}
\label{subsection:formulation:chaser_dynamics}

\begin{figure}
  \centering
  \includegraphics[width=0.6\textwidth]{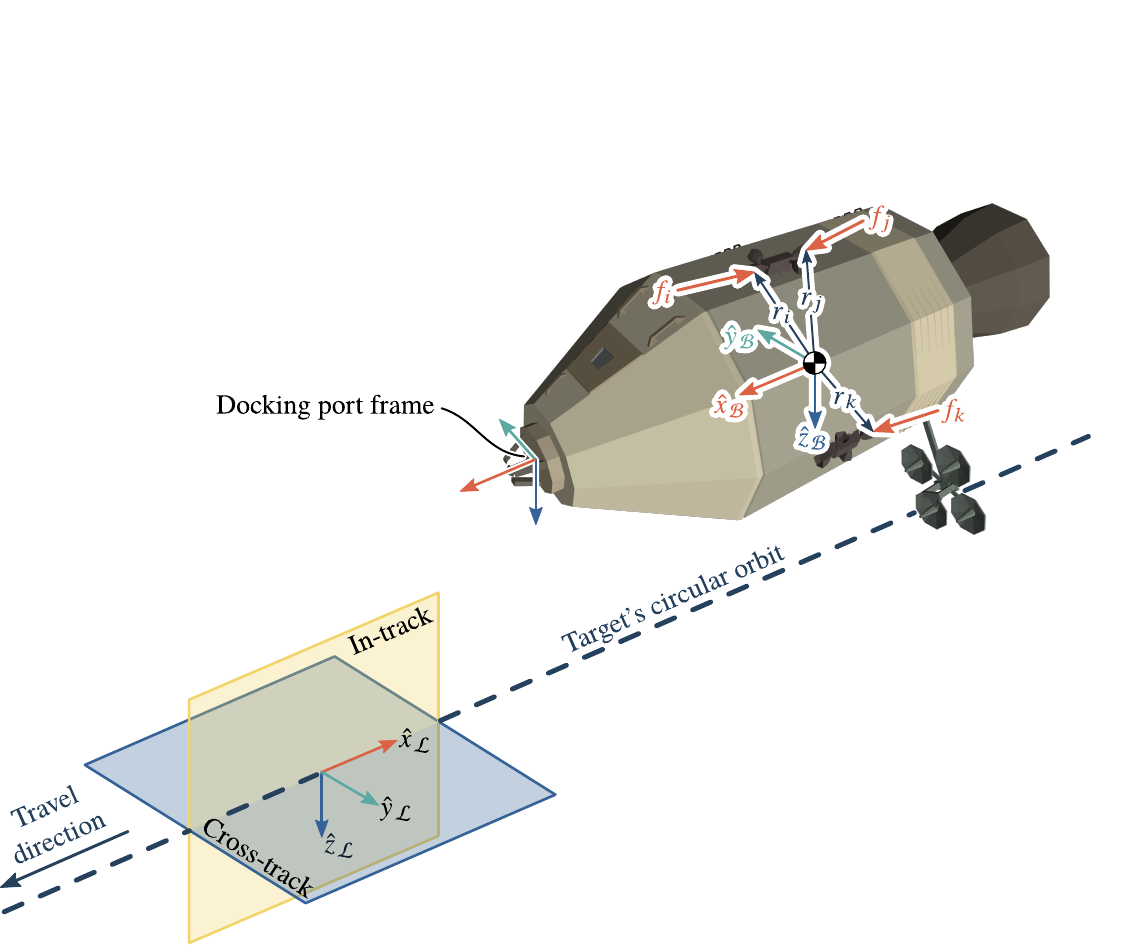}
  \caption{The rendezvous dynamics are written in a Local\dash Vertical
    Local\dash Horizontal frame affixed to the target spacecraft center of
    mass.}
  \label{fig:dynamics_diagram}
\end{figure}

We begin by writing down the equations of motion for the chaser spacecraft. It
is assumed that the chaser is a 6-DOF rigid body vehicle with constant
mass. The latter assumption is accurate for our ultimate numerical application
to the Apollo Transposition and Docking maneuver, whose fuel mass allocation is
$32~\si{\kilo\gram}$, corresponding to about $0.1\%$ of the total Apollo
Command and Service Module (CSM) vehicle mass \cite{NASA_apollo_mass}.

The general setup is illustrated in \figref{dynamics_diagram}. First, a Local\dash
Vertical Local\dash Horizontal (LVLH) frame is placed at the target's center of
mass (COM). Assuming that the target is in a circular orbit, and because
separation distances during the final stages of rendezvous are relatively
small, we can write the translation dynamics in this frame according to the
Clohessy-Wiltshire-Hill equations \cite{CurtisBook}. For the attitude dynamics,
a body frame is affixed to the chaser's COM. Apart from the non\dash inertial
forces of the relative motion dynamics in the LVLH frame, the only forces
acting on the chaser are the ones generated by its system of reaction control
system (RCS) thrusters. As shown in \figref{dynamics_diagram}, the force produced
by each thruster is defined by its point of application $r_i$ and its vector
$f_i$, both of which are expressed in the $\frame{\frB}$ frame. Altogether, the
6-DOF equations of motion of the chaser in the LVLH frame are written as
follows:
\begin{subequations}
  \label{eq:eom}
  \begin{align}
    \label{eq:eom_position}
    \dot p(t) &= v(t), \\
    \label{eq:eom_velocity}
    \dot v(t) &= \frac 1m\sum_{i=1}^\nrcs
                \fade{q(t)\otimes f_i\otimes q(t)\conj}+
                a_{\mrm{LVLH}}\pare[big]{p(t), v(t)}, \\
    \label{eq:eom_rotation}
    \dot q(t) &= \frac 12 q(t)\otimes\omega(t), \\
    \label{eq:eom_rotation_rate}
    \dot \omega(t) &= J\inv\brak[bigg]{
                     \sum_{i=1}^\nrcs
                     \fade{r_i\times f_i(t)}-
                     \omega(t)\times\pare[big]{J\omega(t)}},
  \end{align}
\end{subequations}
where the acceleration due to relative motion is given by:
\begin{equation}
  \label{eq:eom_lvlh_relative_motion}
  a_{\mrm{LVLH}}\pare[big]{p, v} =
  \pare[big]{-2\norb\ez[\frL]\T v}\ex[\frL]+
  \pare[big]{-\norb^2\ey[\frL]\T r}\ey[\frL]+
  \pare[big]{3\norb^2\ez[\frL]\T r+2\norb\ex[\frL]\T v}\ez[\frL],
\end{equation}
where $\norb\in\reals$ is the orbital mean motion. The translation dynamics are
encoded by $p\in\reals^3$ and $v\in\reals^3$, which are LVLH frame vectors
denoting the position and velocity of $\frame{\frB}$ with respect to
$\frame{\frL}$. The attitude dynamics are encoded by a quaternion
$q\in\quaternion$ and an angular velocity $\omega\in\reals^3$. We use the
Hamilton quaternion convention and represent $q$ as a four\dash element vector
\cite{Sola2017}. The quaternion thus represents a frame transformation from
$\frame{\frB}$ to $\frame{\frL}$, or (equivalently) the rotation of a vector in
the $\frame{\frL}$ frame. The $\omega$ vector corresponds to the angular
velocity of $\frame{\frB}$ with respect to $\frame{\frL}$, expressed as a
vector in the $\frame{\frB}$ frame. Altogether, the vehicle state is encoded by
$x=\brak[big]{p; v; q; \omega}\in\reals^{13}$.

\defvar[hide]{$r_i$}{body frame position of the \th{$i$} thruster, \si{\meter}}
\defvar[hide]{$\hat f_i$}{body frame thrust direction of the \th{$i$} thruster}
\defvar[hide]{$\frame{\frL}$}{LVLH frame centered at the target COM}
\defvar[hide]{$\frame{\frB}$}{body frame centered at the chaser COM}
\defvar[hide]{$\norb$}{orbital mean motion, \si{\radian\per\second}}
\defvar[hide]{$m$}{chaser mass, \si{\kilo\gram}}
\defvar[hide]{$J$}{chaser inertia matrix, \si{\kilo\gram\meter\squared}}
\defvar[hide]{$p$}{chaser position, \si{\meter}}
\defvar[hide]{$v$}{chaser velocity, \si{\meter\per\second}}
\defvar[hide]{$q$}{chaser quaternion attitude}
\defvar[hide]{$\omega$}{chaser angular velocity, \si{\radian\per\second}}

\subsection{Impulsive Thrust Model}
\label{subsection:formulation:impulsive_thrust}

As mentioned in the previous section, the chaser is controlled by a system of
\defvar{$\nrcs$}{number of RCS thrusters} RCS thrusters. In accordance with our
ultimate application to the Apollo CSM spacecraft, we assume that each thruster
is able to deliver a constant thrust for a variable duration of time
\cite{NASA_csm_aoh,NASA_apollo_news,NASA_aoh_vol1}. This is known as pulse\dash
width modulation (PWM).

\begin{figure}
  \centering
  \includegraphics{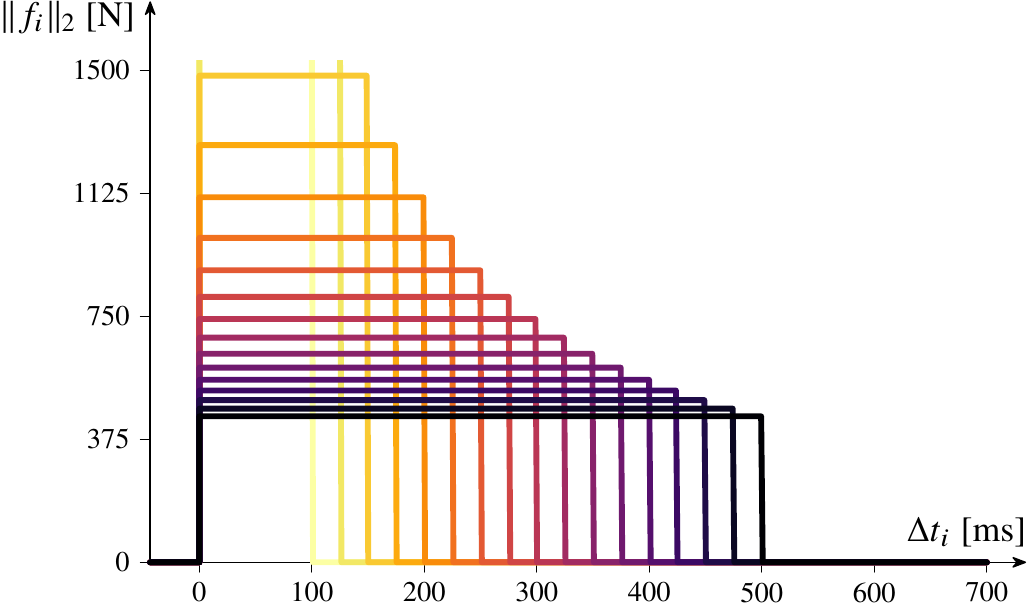}
  \caption{Rectangular thrust pulse duration can be decreased while maintaining
    a constant net impulse (i.e., the area under the curve) by increasing the
    corresponding thrust level.}
  \label{fig:rectangular_pulse_shrinking}
\end{figure}

Let us temporarily focus the discussion on the force produced by the \th{$i$}
thruster. Let $\Frcs$ denote the constant thrust level generated when the
thruster is active (i.e., ``firing''), and let $\Dt[i]$ be the firing or pulse
duration. If the thruster fires for a very short duration relative to the
bandwidth of the chaser's dynamics, then we can approximate the state as being
constant over the firing interval. We can furthermore shrink the firing
interval to zero, as long as we increase the thrust level to maintain a
constant net impulse that is imparted on the chaser. This is illustrated in
\figref{rectangular_pulse_shrinking}, where an original
$500~\si{\milli\second}$ rectangular pulse is reduced down to
$100~\si{\milli\second}$. In the limit as $\Dt[i]$ is reduced to zero, the
thrust signal becomes a scaled Dirac delta function:
\begin{equation}
  \label{eq:dirac_impulse_simple}
  f_i(t) = \Dt[i]\Frcs\delta(t).
\end{equation}

This model is an accurate enough approximation for generating long duration
trajectories with relatively few intermittent control interventions. By
neglecting state variation over the firing duration, the model furthermore has
a significant computational advantage when it comes to linearizing,
discretizing, and simulating the dynamics for the solution process in
\sref{scp}. We emphasize, however, that \eqref{eq:dirac_impulse_simple} is a model
which we use for computation alone. In the physical world, we still expect the
thrusters to fire for a finite duration and at their design (finite) thrust
level.

The discussion so far has centered around a single pulse that occurs at
$t=0~\si{\second}$. We now generalize this model to the trajectory generation
context. Begin by fixing a control interval $\tc>0$ that corresponds to the
``silent'' time interval between thruster firings. Furthermore, let $\Nc$ be
the total number of control opportunities during the trajectory. This means
that the trajectory lasts for $\Nc\tc$ seconds. Note that no firing occurs at
the final time instant, since that would lead to undesirable control at the
moment of docking. Thus, a thruster can be activated only at the time instances
$(k-1)\tc$ where $k=1,2,\dots,\Nc$. To keep the notation short, we define
$\minusone{k}\equiv k-1$ for any general index $k$. Thus, the thrust signal for
the \th{$i$} thruster can be formally written as:
\begin{equation}
  \label{eq:thrust_impulse_signal}
  f_i(t) = \sum_{k=1}^{\Nc}\Dt[ik]\Frcs\delta\pare[big]{t-\minusone{k}\tc}\fdir[i],
\end{equation}
where $\Dt[ik]$ is the pulse duration for the \th{$i$} thruster at the \th{$k$}
control opportunity, and $\fdir[i]$ is the thrust direction unit vector in the
$\frame{\frB}$ frame. Due to delays in on-board electronics and residual
propellant flow downstream of the injector valves
\cite[pp.~2.5-16~to~2.5-18]{NASA_aoh_vol1}, the pulse duration is lower bounded
such that $\Dt[ik]\ge\Dt[\min]$. This is known as a minimum impulse\dash bit
(MIB) constraint. Other propulsion and RCS parameters, such as engine service
life and damage to engine materials, impose an upper bound
$\Dt[ik]\le\Dt[\max]$. As a result, the pulse duration must satisfy the
following nonconvex constraint:
\begin{equation}
  \label{eq:pulse_duration_constraint}
  \Dt[ik]\in\{0\}\union[\Dt[\min], \Dt[\max]].
\end{equation}

\defvar[hide]{$\minusone{k}$}{$k-1$ for any general index $k$}

\begin{figure}
  \centering
  \includegraphics{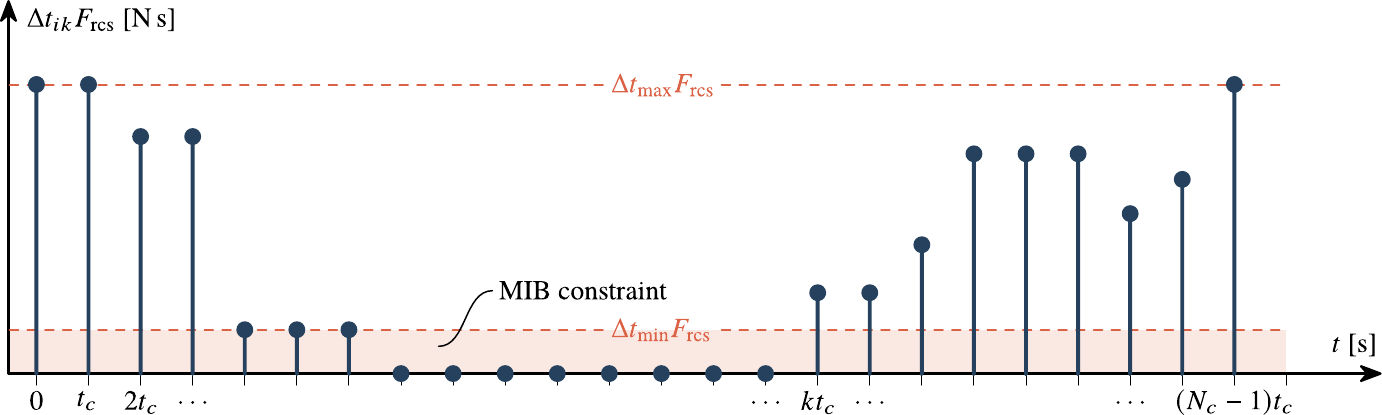}
  \caption{Example of a control history that is compatible with the impulsive
    thrust model \eqref{eq:thrust_impulse_signal} and the pulse duration constraint
    \eqref{eq:pulse_duration_constraint}.}
  \label{fig:pulse_history_illustration}
\end{figure}

\figref{pulse_history_illustration} illustrates a typical control history that we
can expect from the model \eqref{eq:thrust_impulse_signal} subject to the
constraint \eqref{eq:pulse_duration_constraint}. The salient feature of this
control history is that the thruster is either silent, or firing with a minimum
impulse. In particular, no impulse occurs in the MIB keep\dash out zone between
$0$ and $\Dt[\min]\Frcs$~\si{\newton\second}. This region represents impulses
which the RCS system cannot reproduce.

\defvar[hide]{$\tc$}{time interval between thruster firing, \si{\second}}
\defvar[hide]{$\Frcs$}{constant thrust level generated by a thruster, \si{\newton}}
\defvar[hide]{$\fdir[i]$}{thrust direction vector for the \th{$i$} thruster}
\defvar[hide]{$\Dt[ik]$}{pulse duration of \th{$i$} thruster at $k$-th control
  interval, \si{\second}}
\defvar[hide]{$\Dt[\min]$}{minimum pulse duration, \si{\second}}
\defvar[hide]{$\Dt[\max]$}{maximum pulse duration, \si{\second}}
\defvar[hide]{$\Nc$}{number of control opportunities over the trajectory}

\subsection{Plume Impingement Constraint}
\label{subsection:formulation:plume_impingement}

A plume impingement constraint prevents the RCS thrusters from firing and
potentially damaging the target spacecraft. Naturally, this constraint is only
required once the chaser is close enough to the target. Let $\Ifr$ denote the
indices of forward-facing thrusters that are physically pointed along the
$+\ex[\frB]$ axis in \figref{dynamics_diagram}. Due to the physics of rendezvous
and the approach cone constraint of the next section, it is reasonable to
assume that large-angle maneuvering is finished by the time the spacecraft is
close to the target. Thus, when the plume impingement constraint is relevant,
the chaser is approximately facing the target. This yields a simple plume
impingement heuristic: shut off the $\Ifr$ thrusters when the chaser is inside
a so\dash called plume impingement sphere of radius $\rplume$ centered at the
target. This can be formally stated as the following implication:
\begin{equation}
  \label{eq:plume_impingement_constraint}
  \norm[2]{p(\minusone{k}\tc)}\le\rplume~\implies~
  \Dt[ik]=0~\mrm[t]{for all}~i\in\Ifr.
\end{equation}

\defvar[hide]{$\Ifr$}{indices of forward-facing thrusters}
\defvar[hide]{$\rplume$}{plume impingement sphere radius, \si{\meter}}

\subsection{Approach Cone Constraint}
\label{subsection:formulation:approach_cone}

\begin{figure}
  \centering
  \includegraphics{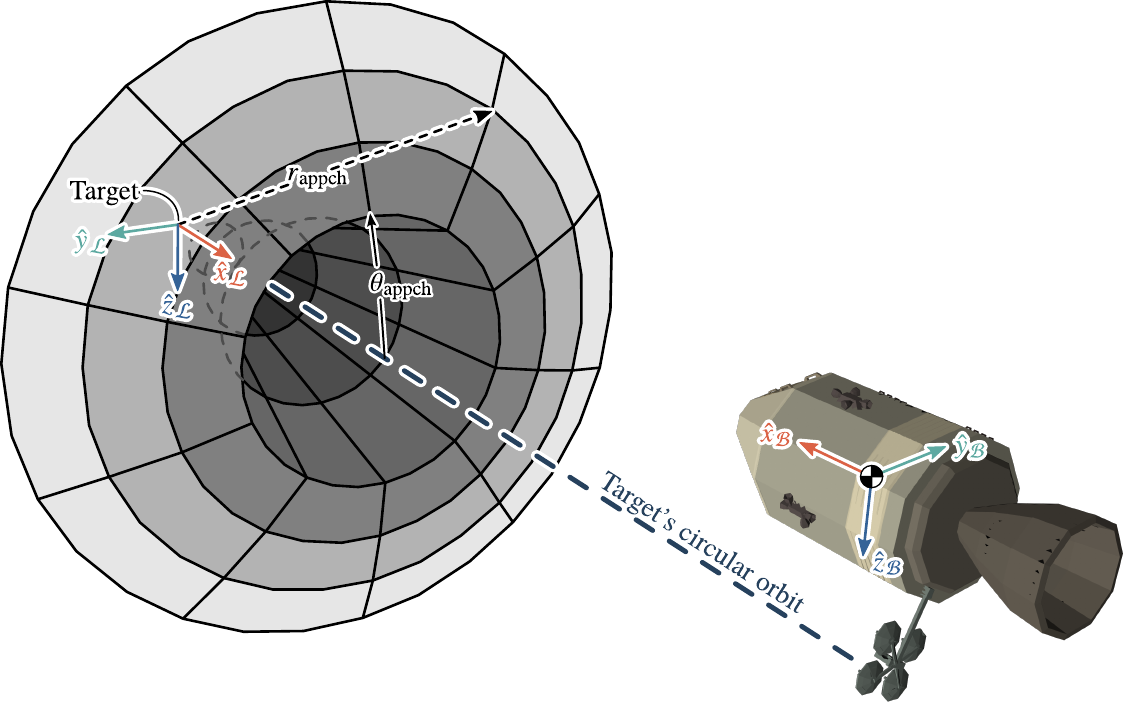}
  \caption{%
    The chaser's position is constrained to lie inside of an approach cone
    when the chaser enters an approach sphere of radius $\rappch$ centered at the
    target (only part of the sphere is shown).}
  \label{fig:approach_cone}
\end{figure}

The approach cone constraint bounds how much the chaser spacecraft can maneuver
once it gets close enough to the target. It has the direct effect of bounding
transverse motion along the $\ey[\frL]$ and $\ez[\frL]$ LVLH axes in
\figref{dynamics_diagram}. In practice, it also bounds all other maneuvering,
including attitude rates, except for translation motion along $-\ex[\frL]$.

\figref{approach_cone} illustrates our implementation of an approach cone. Because
we do not want to restrict the chaser's motion far away from the target, the
constraint only gets applied once the chaser enters a so\dash called approach
sphere of radius $\rappch$. When this condition is satisfied, the chaser's
position is constrained to lie in a cone that emanates from the target along
$+\ex[\frL]$ and has an opening half\dash angle $\angappch$. Formally, the
approach cone constraint can be written as the following implication:
\begin{equation}
  \label{eq:approach_cone_constraint}
  \norm[2]{p(t)}\le\rappch~\implies~
  \ex[\frL]\T p(t)\ge\norm[2]{p(t)}\cos(\angappch).
\end{equation}

\defvar[hide]{$\rappch$}{approach sphere radius, \si{\meter}}
\defvar[hide]{$\angappch$}{approach cone half\dash angle, \si{\radian}}

\subsection{Boundary Conditions}
\label{subsection:formulation:bcs}

We consider the case of terminal rendezvous between two fixed boundary
conditions: some initial chaser state and a terminal ``docked'' state. In
particular, let $x_0=\brak[big]{p_0;v_0;q_0;\omega_0}\in\reals^{13}$ and
$x_f=\brak[big]{p_f;v_f;q_f;\omega_f}\in\reals^{13}$ correspond to the initial
and terminal desired states. The terminal position and attitude are computed
according to the relative geometry of the target and chaser docking ports and
the chaser COM. For simplicity, assume that the target docking port is centered
at the origin of $\frame{\frL}$ and points along $+\ex[\frL]$. Generalizing
this assumption to a non\dash collocated docking port is possible, but does not
represent anything particularly novel or challenging for our algorithm. When
docked, let $\qL\in\quaternion$ denote the chaser docking port's attitude with
respect to the target docking port. As illustrated in \figref{dynamics_diagram},
$\qL$ is a simple yaw around $+\ez[\frL]$ by $180\si{\degree}$. Furthermore,
let $\qDP\in\quaternion$ and $\rDP\in\reals^3$ be the rotation and position of
the chaser docking port relative to $\frame{\frB}$. The terminal position and
attitude are then given by:
\begin{subequations}
  \begin{align}
    \label{eq:terminal_quatenion}
    q_f &= \qL\otimes\qDP\conj, \\
    \label{eq:terminal_position}
    p_f &= -q_f\otimes\rDP\otimes q_f\conj.
  \end{align}
\end{subequations}

\defvar[hide]{$\tf$}{trajectory duration, \si{\second}}
\defvar[hide]{$p_0,v_0,q_0,\omega_0$}{initial boundary conditions}
\defvar[hide]{$p_f,v_f,q_f,\omega_f$}{terminal boundary conditions}
\defvar[hide]{$\tol{p_f},\tol{v_f},\tol{q_f},\tol{\omega_f}$}{%
    terminal boundary condition tolerances}

For a rendezvous trajectory that lasts $\tf$ seconds, the boundary conditions
we impose are:
\begin{equation}
  \label{eq:boundary_conditions}
  x(0) = x_0,\quad x(\tf)+\Dxf = x_f,
\end{equation}
where
$\Dxf=\brak[big]{\Relax{p_f};\Relax{v_f}; \Relax{q_f}; \Relax{\omega_f}}\in\reals^{13}$
relaxes of the terminal boundary condition. This is necessary because the MIB
constraint from \figref{pulse_history_illustration} makes it impossible to
fine\dash tune the trajectory to arbitrary precision. In general, some terminal
error has to occur. As long as this error is small, it will be safely absorbed
by the mechanical design of the docking port. The required tolerances can be
found in the spacecraft's manual. For example, for the Apollo CSM the complete
list is given in \cite[Section~3.8.2.3]{NASA_csm_aoh}. Because it is good
practice to leave a margin of error for feedback controllers, we will constrain
$\Dxf$ to a much smaller value than what the docking mechanism can
tolerate. The following constraints restrict the size of $\Dxf$ to user\dash
specified tolerances:
\begin{subequations}
  \begin{align}
    \label{eq:pf_tolerance}
    \norm[\infty]{\Relax{p_f}}
    &\le \tol{p_f},\quad\ex[\frL]\T\Relax{p_f} = 0, \\
    \label{eq:vf_tolerance}
    \norm[\infty]{\Relax{v_f}} &\le \tol{v_f}, \\
    \label{eq:qf_tolerance}
    q(\tf)\T q_f &\ge \cos\pare{\tol{q_f}/2}, \\
    \label{eq:wf_tolerance}
    \norm[\infty]{\Relax{\omega_f}} &\le \tol{\omega_f}.
  \end{align}
\end{subequations}

The terminal position along $\ex[\frL]$ is made exact since contact along
$\ex[\frL]$ is required for docking. Furthermore, it is always possible to
satisfy by adjusting $\tf$. The terminal attitude is constrained by
\eqref{eq:qf_tolerance} in terms of an error quaternion, and says that the angular
deviation from $q_f$ about any axis must be no larger than angle $\tol{q_f}$.

\subsection{Basic Rendezvous Problem}
\label{subsection:formulation:summary}

Our goal is to compute a fuel\dash optimal rendezvous trajectory, which means
that it is desirable to keep the pulse durations $\Dt[ik]$ as short and as
sparse as possible. An appropriate optimization cost function is simply the sum
of pulse durations for all thrusters and control opportunities:
\begin{equation}
  \label{eq:original_cost_function}
  \Jfuel = \Dt[\max]\inv\sum_{i=1}^{\nrcs}\sum_{k=1}^{\Nc}\Dt[ik],
\end{equation}
where the normalization by $\Dt[\max]$ is useful when
\eqref{eq:original_cost_function} is mixed with other costs for the solution
process in \sref{scp}. Note that \eqref{eq:original_cost_function} is effectively a
one-norm penalty on the pulse durations. This encourages the optimal pulse
history to be sparse, which goes part of the way towards discouraging MIB
constraint violation \cite{BoydConvexBook,Hartley2013}.

We can now summarize the above sections by writing the full rendezvous
optimization problem that has to be solved. We call this the basic rendezvous
problem (BRP). Starting now and throughout the rest of the article, the time
argument will be omitted whenever it does not introduce ambiguity.
\begin{optimization}[
  label={brp},
  variables={x, \Dt, \tf},
  objective={\Jfuel}]
  \optilabel{position}
  &\dot p = v, \\
  &\dot v = \frac 1m\sum_{i=1}^\nrcs
  \fade{q\otimes f_i\otimes q\conj}+a_{\mrm{LVLH}}\pare[big]{p, v}, \\
  &\dot q = \frac 12 q\otimes\omega, \\
  \optilabel{angular_velocity}
  &\dot \omega = J\inv\brak[bigg]{
    \sum_{i=1}^\nrcs
    \fade{r_i\times f_i}-
    \omega\times\pare[big]{J\omega}}, \\
  \optilabel{mib}
  &\Dt[ik]\in\{0\}\union[\Dt[\min], \Dt[\max]], \\
  \optilabel{plume}
  &\norm[2]{p(\minusone{k}\tc)}\le\rplume~\implies~
  \Dt[ik]=0~\mrm[t]{for all}~i\in\Ifr, \\
  \optilabel{approach}
  &\norm[2]{p}\le\rappch~\implies~
  \ex[\frL]\T p\ge\norm[2]{p}\cos(\angappch), \\
  \optilabel{bcs1}
  &x(0) = x_0,\quad x(\tf)+\Dxf = x_f, \\
  &\norm[\infty]{\Relax{p_f}} \le \tol{p_f},~\ex[\frL]\T\Relax{p_f} = 0,~
  \norm[\infty]{\Relax{v_f}} \le \tol{v_f}, \\
  \optilabel{bcs3}
  &q(\tf)\T q_f \ge \cos\pare{\tol{q_f}/2},~
  \norm[\infty]{\Relax{\omega_f}} \le \tol{\omega_f}.
\end{optimization}

The BRP is a continuous\dash time, free\dash final time, nonconvex optimal
control problem. It is not efficiently solvable on a computer for three main
reasons \cite{CSM2021}:
\begin{enumerate}
\item Continuous\dash time problems have an infinite number of DOFs in the
  optimized control signal. However, numerical optimization algorithms are
  restricted to a finite number of DOFs;
\item The problem has nonlinear dynamics, which results in a nonconvex optimal
  control problem. However, numerical algorithms for nonconvex problems require
  expert initial guesses and generally do not converge quickly and reliably
  enough for safety\dash critical applications
  \cite{NocedalBook,BoydConvexBook};
\item The constraints \optieqref{brp}{mib}-\optieqref{brp}{approach} contain discrete
  \textifelse logic. This is traditionally handled by mixed\dash integer
  programming (MIP), which has exponential computational complexity and does
  not scale well to large problems \cite{achterberg2013progress}.
\end{enumerate}

We will begin by resolving the third issue through a homotopy approach in the
next section. The first two issues will then be tackled in \sref{scp}.

\section{Homotopy for Smooth Approximation of Discrete Logic}
\label{section:embedding}

We now consider the problem of computationally efficient modeling the discrete
logic constraints \optieqref{brp}{mib}-\optieqref{brp}{approach}. This model along with
the associated numerical continuation solution method in \sref{scp} are the main
contributions of this paper. We begin in \ssref{embedding:motivation} with a
motivation for why a new approach to handling discrete logic is necessary. Our
homotopy approach is then described in general terms in
\ssref{embedding:description}. Finally,
\ssref{embedding:approach,embedding:impingement,embedding:mib} specialize the
approach to the discrete logic constraints
\optieqref{brp}{mib}-\optieqref{brp}{approach}.

\subsection{Motivation}
\label{subsection:embedding:motivation}

The traditional way of handling discrete logic in an optimization problem is
through the use of binary variables \cite{Schouwenaars2001,Schouwenaars2006}. As
a concrete example, consider the plume impingement constraint
\optieqref{brp}{plume}. Let $\zplume(t):[0,\tf]\to\{0,1\}$ denote a binary variable
trajectory that is also to be optimized. Let $\Mplume$ be a large positive
value that bounds all possible values of $\norm[2]{p(t)}$ that can occur during
a rendezvous trajectory. For example, $\Mplume=10\norm[2]{p_0}$ is a reasonable
choice. The plume impingement constraint \optieqref{brp}{plume} can then be
equivalently written as:
\begin{subequations}
  \label{eq:mip_plume_impingement}
  \begin{align}
    \label{eq:mip_plume_impingement_trigger}
    \zplume(\minusone{k}\tc)\rplume
    &\le \norm[2]{p(\minusone{k}\tc)}
      \le \rplume+\zplume(\minusone{k}\tc)\Mplume, \\
    \label{eq:mip_plume_impingement_constraint}
    0 &\le \Dt[ik] \le \zplume(\minusone{k}\tc)\Dt[\max]~\mrm[t]{for all}~
        i\in\Ifr.
  \end{align}
\end{subequations}

Looking at \eqref{eq:mip_plume_impingement}, $\zplume$ can be interpreted as
follows: the chaser is outside the plume impingement sphere if and only if
$\zplume=1$. When the chaser is inside this sphere, the only feasible choice is
$\zplume=0$, and \eqref{eq:mip_plume_impingement_constraint} shuts off the
forward\dash facing thrusters.

A similar formulation can be used to model the MIB and approach cone
constraints \optieqref{brp}{mib} and \optieqref{brp}{approach}, resulting in a MIP
formulation. Unfortunately, this approach has an issue when it comes to
actually solving \pref{brp}: mixed\dash integer optimization algorithms are
generally too slow for real\dash time applications, are computationally
expensive, and do not scale well to large problem sizes
\cite{achterberg2013progress,MalyutaLCSS}. When compounded by the fact that this
formulation introduces new nonconvex constraints (e.g., the position norm lower
bound in \eqref{eq:mip_plume_impingement_trigger}), it becomes clear that the MIP
approach is not a workable real\dash time solution method for \pref{brp}.

Several methods have been proposed in recent years to replace MIP with a
real\dash time capable approach. On the one hand, recent theoretical results
have demonstrated that a lossless relaxation can solve certain classes of
problems with discrete logic constraints on the control variable
\cite{Harris2021,MalyutaIFAC}. This approach is practical because it requires
solving only a single convex problem. Some versions of the method can handle
restricted forms of nonlinear dynamics and convex constraints
\cite{Blackmore2012,HarrisThesis,CSM2021}. However, the method does not apply to
the full generality of \pref{brp}, which involves more complicated nonlinear
dynamics as well as discrete logic constraints on the state.

A separate family of solution methods has been proposed to handle discrete
logic constraints using sequential convex programming (SCP) \cite{CSM2021}. The
methods define so\dash called state triggered constraints (STCs) that can embed
general discrete logic into a continuous optimization framework
\cite{SzmukThesis,ReynoldsThesis,ARC2021}. Two equivalent forms of STCs have
been proposed, based on a slack variable \cite{Szmuk2020} and based on a
multiplicative coefficient that is motivated by the linear complementarity
problem \cite{Reynolds2020}. STCs have also been extended to handle quite
general logical combinations of \textand and \textor gates
\cite{Szmuk2019a,Szmuk2019b}. In fact, the authors of this paper have applied
STCs to solve a version of \pref{brp}, with the results available in
\cite{MalyutaScitechDocking}. In the latter work it was observed that STCs run
into an issue called \textit{locking} for the MIB constraint \optieqref{brp}{mib}
\cite[Definition~1]{MalyutaScitechDocking}. As described in \sref{scp}, SCP works by
iteratively refining an approximate solution of \pref{brp}. In brief terms, locking
means that once the algorithm chooses $\Dt[ik]=0$ at a particular iteration, it
is unable to change the value to $\Dt[ik]\in[\Dt[\min],\Dt[\max]]$ at later
iterations. The effect is that the algorithm is susceptible to getting into a
``corner'' where it is unable to use thrusters if they become needed at later
refinements of the rendezvous trajectory. The consequence is failure to
generate a feasible trajectory. There is currently no known remedy for
constraints that exhibit locking in the STC formulations of
\cite{Szmuk2020,Reynolds2020}.

For reasons that are well documented in past literature, we view SCP as one of
the most effective frameworks for the real\dash time solution of nonconvex
trajectory generation problems \cite{CSM2021,ARC2021,Reynolds2020Real}. Thus,
our primary motivation is to devise a new general method that is free from
locking and that can embed discrete logic into an SCP-based continuous
optimization framework.

\subsection{The Homotopy Algorithm}
\label{subsection:embedding:description}

We now develop a homotopy\dash based method to systematically handle
\textifelse discrete logic constraints of the following form:
\begin{subequations}
  \label{eq:discrete_logic_couple}
  \begin{align}
    \label{eq:if_implication}
    \mrm[t]{``If''}\quad \Lif(z)\definedas
    \bigwedge_{i=1}^{\ng}\pare[bigg]{\predicate[i](z)\le 0}~
    &\implies~\fL(z)\le 0, \\
    \label{eq:else_implication}
    \mrm[t]{``Else''}\quad \Rif(z)\definedas
    \bigvee_{i=1}^{\ng}\pare[bigg]{\predicate[i](z)> 0}~
    &\implies \fR(z)\le 0,
  \end{align}
\end{subequations}
where $z\in\reals^{\nz}$ is a generic placeholder for one or several
optimization variables. The functions $\predicate[i]:\reals^{\nz}\to\reals$ are
called predicates, and the functions $\fL:\reals^{\nz}\to\reals^{\nf}$ and
$\fR:\reals^{\nz}\to\reals^{\nf}$ are implication constraints to be enforced
when the corresponding expression's left\dash hand side is \textcode{true}. For
\eqref{eq:if_implication} this is a combination of \textand gates, whereas for
\eqref{eq:else_implication} it is a combination of \textor gates with the predicate
inequalities reversed. We may thus see \eqref{eq:discrete_logic_couple} in the
following light: enforce $\fL\le 0$ when \textit{all} the predicates are
nonpositive, or enforce $\fR\le 0$ when \textit{any} predicate is positive.

\begin{figure}
  \centering
  \includegraphics[width=1\textwidth]{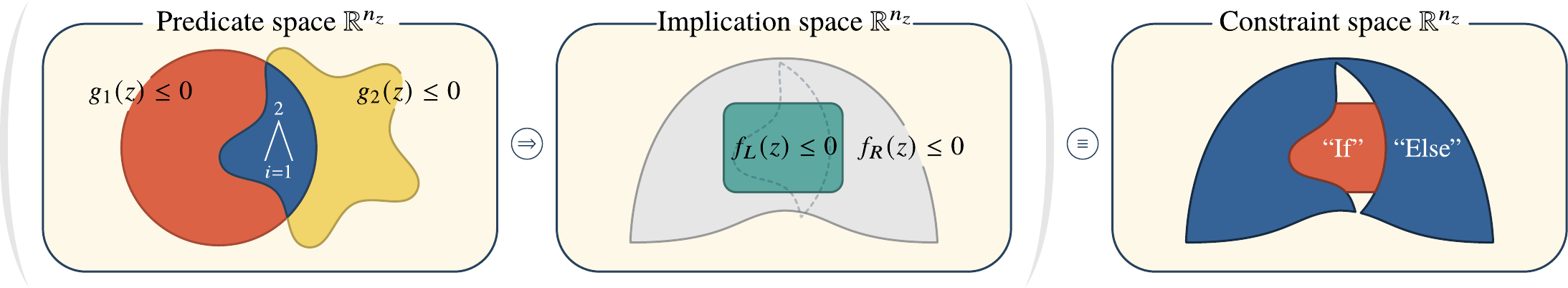}
  \caption{%
    Pictorial representation of the \textifelse discrete logic constraint
    \eqref{eq:discrete_logic_couple}.
  }
  \label{fig:discrete_logic_illustration}
\end{figure}

One can show using De Morgan's theorem that $\Rif(z)=\neg \Lif(z)$. As a
result, the implications in \eqref{eq:discrete_logic_couple} indeed form an
\textifelse pair in the sense that exactly one of $\fL$ and $\fR$ is enforced
at any given instant. The situation is illustrated in
\figref{discrete_logic_illustration}. In the predicate space, the functions
$\predicate[i]$ form sublevel sets of nonpositive values. In the implication
space, the constraint functions $\fL$ and $\fR$ also form sublevel sets of
nonpositive values. Note that these sets can generally be disjoint. The overall
\textif constraint is obtained by intersecting the sublevel set of $\fL$ with
the sublevel set of the \textand combination. Conversely, the overall \textelse
constraint is obtained by intersecting the sublevel set of $\fR$ with the
complement of the sublevel set for the \textand combination.

By using the value 1 to denote \texttrue and 0 to denote \textfalse, we have
the complementarity relationship $\Rif(z)=1-\Lif(z)$. Using this property,
\eqref{eq:discrete_logic_couple} can be stated in the following equivalent ways:
\begin{subequations}
  \label{eq:discrete_logic_exact}
  \begin{align}
    \label{eq:discrete_logic_exact_and}
    \Lif(z)\fL(z)+\brak[big]{1-\Lif(z)}\fR(z) &\le 0, \\
    \label{eq:discrete_logic_exact_or}
    \brak[big]{1-\Rif(z)}\fL(z)+\Rif(z)\fR(z) &\le 0.
  \end{align}
\end{subequations}

Because \eqref{eq:discrete_logic_exact} involves discrete elements (i.e., the
\textand and \textor gates), it cannot be readily included in a continuous
optimization problem. As mentioned in the previous section, STCs are one
possible way to circumvent the issue, however they exhibit locking in the
particular case of the MIB constraint \optieqref{brp}{mib}. An alternative approach
is to replace either $\Lif$ or $\Rif$ by a smooth approximation, and to apply a
numerical continuation scheme to iteratively improve the approximation until
some arbitrary precision \cite{NocedalBook}. We take this latter approach, and
begin with a brief description of two existing methods.

\subsubsection{Existing Homotopy Methods}
\label{subsubsection:embedding:description:existing}

Homotopy is the core idea behind the recent relaxed autonomous switched hybrid
system (RASHS) and composite smooth control (CSC) algorithms
\cite{saranathan2018relaxed,taheri2020novel,arya2021composite}. Both algorithms
model the constraint \eqref{eq:discrete_logic_exact_and} by approximating the
\textand combination with a sigmoid function. To this end, let
$\sigma_{\sharp}(w):\reals\to\reals$ represent a sigmoid function which
approaches one for negative arguments and zero for positive arguments. The
transition point occurs at $w=0$ and the homotopy parameter $\sharp>0$ (also
known as a sharpness parameter) regulates how quickly the transition
happens. As $\sharp$ increases, $\sigma_{\sharp}$ approaches a ``step down''
function. This allows RASHS and CSC to model $\Lif$ as follows:
\begin{equation}
  \label{eq:discrete_logic_and_approx}
  \Lif(z)\approx\smooth{\Lif}(z)\definedas
  \prod_{i=1}^{\ng} \sigma_{\sharp}\pare[big]{\predicate[i](z)}.
\end{equation}

By replacing $\Lif$ with $\smooth{\Lif}$ in \eqref{eq:discrete_logic_exact_and},
the RASHS and CSC methods can model discrete logic in a smooth way that is
conducive for continuous optimization. By using numerical continuation to
progressively increase $\sharp$, the methods can enforce the discrete logic
constraint \eqref{eq:discrete_logic_exact_and} with arbitrary accuracy.

\subsubsection{Proposed Homotopy Method}
\label{subsubsection:embedding:description:our}

\begin{figure}
  \centering
  \includegraphics{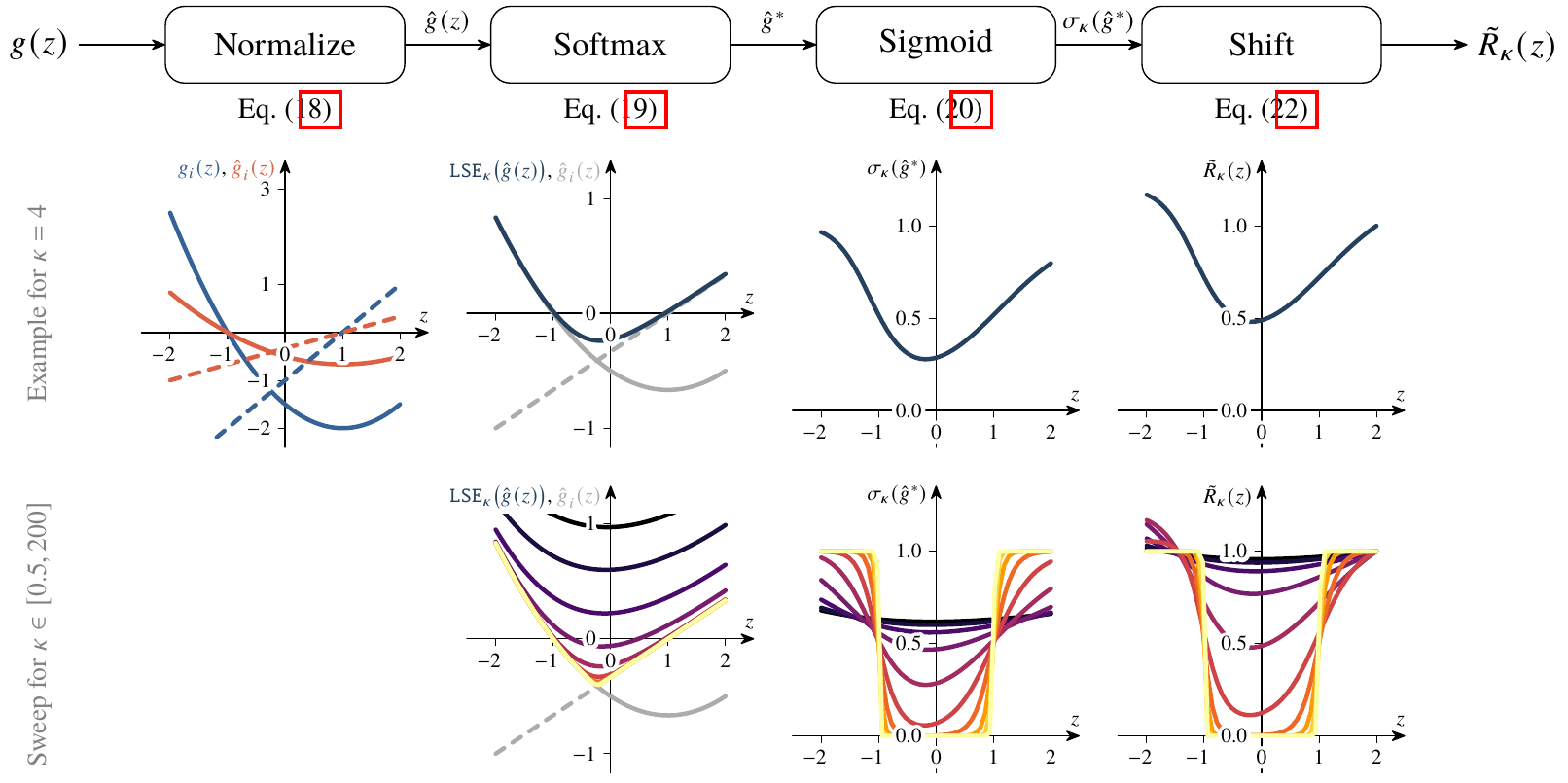}
  \caption{Illustration for a simple set of two predicates ({\color{Blue}blue}
    curves in the second row, first column) of the four stages comprising our
    smooth approximation of the \textor combination in
    \eqref{eq:else_implication}.}
  \label{fig:smooth_approximation_stages}
\end{figure}

Our method for imposing \eqref{eq:discrete_logic_couple} is centered around a
smooth approximation of the alternative constraint
\eqref{eq:discrete_logic_exact_or} using a multinomial logit function
\cite{Hastie2009}. We thus view our approach as a ``dual'' formulation to RASHS
and CSC: instead of modeling the \textand combination of \eqref{eq:if_implication},
we model its complement \eqref{eq:else_implication}. A noteworthy benefit of this
approach is the ability to model \textor combinations, whereas RASHS and CSC
both appear to be compatible only with \textand logic. Our method is therefore
an extension of the ideas in RASHS and CSC. Although we do not develop the full
theory here, our method together with \eqref{eq:discrete_logic_and_approx} can
model arbitrary combinations of \textand and \textor logic. This extends smooth
modeling of discrete logic to its full generality.

We break down the smooth approximation of $\Rif$ into four computational
``stages''. The reader may follow along with the help of the illustration in
\figref{smooth_approximation_stages}. Begin with the raw data, which are the
individual predicate values $\predicate[i](z)$. For convenience, let
$\predicate(z):\reals^{\nz}\to\reals^{\ng}$ be the concatenated vector of
predicates. The first stage is to normalize $g(z)$ by the expected maximum
value of the predicates:
\begin{equation}
  \label{eq:predicate_normalization}
  \predicate[\max]\definedas\max_{z} \norm[\infty]{\predicate(z)},
\end{equation}
where $z$ is understood to be taken from the set of all reasonable values for
\pref{brp}. We can then define a normalized predicate vector:
\begin{equation}
  \label{eq:normalized_predicate}
  \hat \predicate(z) \definedas \predicate[\max]\inv\predicate(z).
\end{equation}

Normalization ensures that $\hat g(z)$ takes values in a $[-1, 1]^{\ng}$
hypercube. This helps to standardize the parameter choices for the numerical
continuation solution method, which we will describe in \sref{scp}. The second
stage is to pick out the maximum predicate value. Because we want a smooth
approximation, we find an approximate maximum using the log\dash sum\dash exp
function, also known as a softmax. For a given homotopy parameter $\sharp>0$,
the softmax function $\softmax:\reals^{\ng}\to\reals$ is defined by:
\begin{equation}
  \label{eq:softmax}
  \softmax\pare[big]{\hat\predicate(z)} \definedas \sharp\inv\log\pare[biggg]{
    \sum_{i=1}^{\ng}\exp{\sharp \hat\predicate_i(z)}}.
\end{equation}
Let us denote the resulting value by
$\hat\predicate^*\equiv\softmax\pare[big]{\hat\predicate(z)}$. As $\sharp$
grows, this value approaches the true $\max_i\hat\predicate_i(z)$. In the third
stage, the value is passed to a sigmoid function which maps it to the $[0,1]$
interval. This function approaches zero for negative arguments and one for
positive arguments. We define it as follows:
\begin{equation}
  \label{eq:sigmoid}
  \sigma_{\sharp}(\hat\predicate^*) \definedas
  1-\brak[bigg]{1+\exp{\sharp\hat\predicate^*}}\inv.
\end{equation}

Note that by subtituting \eqref{eq:softmax} into \eqref{eq:sigmoid}, we obtain the
familiar multinomial logit function \cite{Hastie2009}:
\begin{equation}
  \label{eq:logistic_sigmoid}
  \sigma_{\sharp}(\hat\predicate^*) =
  1-\brak[biggg]{1+\sum_{i=1}^{\ng}\exp{\sharp\hat\predicate_i}}\inv.
\end{equation}

For this reason, we call our approach \textit{multinomial logit
  smoothing}. When $\sharp$ is large and the time comes to computing the
derivatives of \eqref{eq:logistic_sigmoid} for the solution process in \sref{scp}, we
have noted that there are important numerical stability advantages to breaking
the logistic function into separate steps \eqref{eq:softmax} into
\eqref{eq:sigmoid}. This is why we keep the second and third stages separate.

The fourth and last stage of approximating $\Rif$ is to vertically shift the
sigmoid function so that it matches its exact value at some specified predicate
value $\predmatch\in\reals^{\ng}$, where we require at least one element to be
positive (such that $\Rif(z)>0$). We typically choose
$\predmatch=\predicate(z^*)$ where $z^*=\argmax_z\norm[\infty]{\predicate(z)}$
from \eqref{eq:predicate_normalization}. Shifting carries the benefit of not
over\dash restricting the solution variables early in the solution process,
when $\sharp$ is small and $\sigma_{\sharp}\approx \ng/(\ng+1)$. The latter
effect is visible in the bottom row, third column of
\figref{smooth_approximation_stages}. Ultimately, the smooth approximation of
$\Rif$ is defined as follows, and is the direct counterpart of the RASHS and
CSC model \eqref{eq:discrete_logic_and_approx}:
\begin{equation}
  \label{eq:discrete_logic_or_approx}
  \Rif(z)\approx\smooth{\Rif}(z)\definedas
  \sigma_{\sharp}(\hat\predicate^*)+\pare[big]{1-\sigma_{\sharp}(\predmatch)}.
\end{equation}

The discrete logic constraint \eqref{eq:discrete_logic_couple} can then be written
as the following smooth approximation, which is obtained by substituting $\Rif$
in \eqref{eq:discrete_logic_exact_or} with $\smooth{\Rif}$ from
\eqref{eq:discrete_logic_or_approx}:
\begin{equation}
  \label{eq:smooth_discrete_logic}
  \brak[big]{1-\smooth{\Rif}(z)}\fL(z)+\smooth{\Rif}(z)\fR(z) \le 0.
\end{equation}

In the following sections, we will show how to use \eqref{eq:smooth_discrete_logic}
to model the discrete logic constraints
\optieqref{brp}{mib}-\optieqref{brp}{approach}. For the sake of comparison, the RASHS
and CSC smooth approximations \eqref{eq:discrete_logic_and_approx} are given by
\cite{saranathan2018relaxed,taheri2020novel}:
\begin{subequations}
  \label{eq:rashs_csc_smooth_and}
  \begin{align}
    \label{eq:rashs_and}
    \smooth{\Lif}^{\mrm{RASHS}}(z)
    &= \prod_{i=1}^{\ng} \pare[big]{1+\exp{\sharp \hat\predicate_i(z)}}\inv, \\
    \label{eq:csc_and}
    \smooth{\Lif}^{\mrm{CSC}}(z)
    &= \prod_{i=1}^{\ng} \frac 12\pare[bigg]{1-
      \tanh\pare[big]{\sharp \hat\predicate_i(z)}}.
  \end{align}
\end{subequations}

\begin{figure}
  \centering
  \begin{subfigure}[t]{0.33\textwidth}
    \includegraphics[page=1]{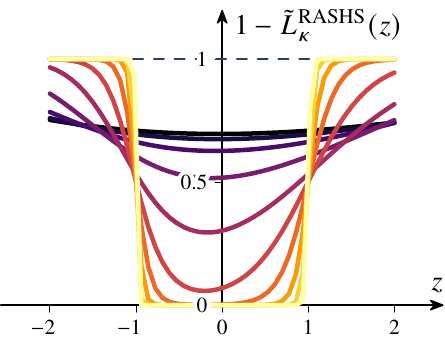}
    \caption{The RASHS method \cite{saranathan2018relaxed}.}
    \label{fig:rashs_smooth_example}
  \end{subfigure}%
  \hfill%
  \begin{subfigure}[t]{0.33\textwidth}
    \includegraphics[page=2]{sigmoid_comparison.pdf}
    \caption{The CSC method \cite{taheri2020novel,arya2021composite}.}
    \label{fig:csc_smooth_example}
  \end{subfigure}%
  \hfill%
  \begin{subfigure}[t]{0.33\textwidth}
    \includegraphics[page=3]{sigmoid_comparison.pdf}
    \caption{Our proposed method \eqref{eq:discrete_logic_or_approx}.}
    \label{fig:us_smooth_example}
  \end{subfigure}
  \caption{%
    Comparison of smoothed discrete \textor logic \eqref{eq:else_implication}
    obtained using three smoothing methods. Faint dashed lines in
    (\protect\subref{fig:us_smooth_example}) show the multinomial logit function
    \eqref{eq:logistic_sigmoid} without shifting, which is very similar to RASHS in
    (\protect\subref{fig:rashs_smooth_example}).}
  \label{fig:smooth_comparison}
\end{figure}

\figref{smooth_comparison} compares the smooth logic \eqref{eq:rashs_csc_smooth_and}
with our approach \eqref{eq:discrete_logic_or_approx}. Without the shifting
operation in \eqref{eq:discrete_logic_or_approx}, all three methods are remarkably
similar. Multinomial logit smoothing without shifting is most similar to RASHS:
the two methods are identical for $\ng=1$, and slightly different for
$\ng>1$. Thus, shifting is the critical difference in our method. As we shall
see below, it is most critical for constraints like the MIB \optieqref{brp}{mib},
where it is important that $\smooth{\Rif}(z)\approx 1$ for small $\sharp$ (this
effectively removes the MIB constraint from the early solution algorithm
iterations in \sref{scp}).

\subsection{Modeling the Approach Cone}
\label{subsection:embedding:approach}

We begin by modeling the approach cone constraint \optieqref{brp}{approach} in the
framework of \eqref{eq:discrete_logic_couple} and its smooth approximation
\eqref{eq:smooth_discrete_logic}. Comparing \optieqref{brp}{approach} with
\eqref{eq:if_implication}, we have $\ng=1$, $z=p$, and the predicate:
\begin{equation}
  \label{eq:approach_cone_predicate}
  \predicate[1](p) = p\T p-\rappch^2,
\end{equation}
where we use the two\dash norm squared to conveniently make the predicate
everywhere smooth. This predicate is then used in
\eqref{eq:discrete_logic_or_approx} to form $\smoothRappch$, the smooth \textor
approximation for the approach cone predicate. The \textif implication can be
written as:
\begin{equation}
  \label{eq:approach_cone_if}
  \fL(p) = \cos(\angappch)-\ex[\frL]\T p\norm[2]{p}\inv.
\end{equation}

When the chaser is outside of the approach sphere, we wish to allow the
chaser's trajectory to assume any approach angle. By the Cauchy\dash Schwarz
inequality, this can be expressed as the inequality
$\ex[\frL]\T p\ge -\norm[2]{p}$. As a result, the \textelse implicationx can be
written as:
\begin{equation}
  \label{eq:approach_cone_else}
  \fR(p) = -1-\ex[\frL]\T p\norm[2]{p}\inv.
\end{equation}

We can now use \eqref{eq:approach_cone_if} and \eqref{eq:approach_cone_else} directly
in \eqref{eq:smooth_discrete_logic}, which yields a smooth approximation of the
approach cone constraint:
\begin{equation}
  \label{eq:smooth_approach_cone}
  \cos(\angappch)-\pare[big]{1+\cos(\angappch)}\smoothRappch(p)-
  \ex[\frL]\T p\norm[2]{p}\inv\le 0.
\end{equation}

\subsection{Modeling Plume Impingement}
\label{subsection:embedding:impingement}

The plume impingement constraint \optieqref{brp}{plume} is modeled in a very
similar way. Recall that the rendezvous trajectory has $\Nc$ control
opportunities and the chaser has $\Ifr$ forward\dash facing thrusters. Let us
focus on the \th{$k$} control opportunity for thruster $i\in\Ifr$. Comparing
\optieqref{brp}{plume} with \eqref{eq:if_implication}, we have $\ng=1$ and
$z=\brak[big]{p(\minusone{k}\tc); \Dt[ik]}$. The predicate takes after
\eqref{eq:approach_cone_predicate}:
\begin{equation}
  \label{eq:plume_predicate}
  \predicate[1]\pare[big]{p(\minusone{k}\tc)} =
  p(\minusone{k}\tc)\T p(\minusone{k}\tc)-\rplume^2,
\end{equation}

This predicate is then used in \eqref{eq:discrete_logic_or_approx} to form
$\smoothRplume$, the smooth \textor approximation for the plume impingement
predicate. The \textif implication for plume impingement is an equality
constraint, whereas our standard formulation \eqref{eq:discrete_logic_couple}
requires an inequality. To reconcile the two situations, one possible approach
is to leverage \optieqref{brp}{mib} and to realize that
$\Dt[ik]\in [0,\Dt[\max]]$. Thus, we can impose the constraint:
\begin{equation}
  \label{eq:pulse_bounds}
  0\le \Dt[ik]\le \Dt[\max],
\end{equation}
and we can write the following \textif implication:
\begin{equation}
  \label{eq:plume_if}
  \fL(\Dt[ik]) = \Dt[ik].
\end{equation}

Equation \eqref{eq:plume_if} together with \eqref{eq:pulse_bounds} enforce
$0\le \Dt[ik]\le 0$ when the predicate \eqref{eq:plume_predicate} is \texttrue,
which is equivalent to \optieqref{brp}{plume}. When the chaser is outside of the
plume impingement sphere, the forward\dash facing thrusters are free to
fire. We can express this as the following \textelse implication:
\begin{equation}
  \label{eq:plume_else}
  \fR(\Dt[ik]) = \Dt[ik]-\Dt[\max].
\end{equation}

Equations \eqref{eq:plume_if} and \eqref{eq:plume_else} can now be substituted into
\eqref{eq:smooth_discrete_logic}, yielding a smooth approximation of the plume
impingement constraint:
\begin{equation}
  \label{eq:smooth_approach_cone}
  \Dt[ik] \le \smoothRplume\pare[big]{p(\minusone{k}\tc)}\Dt[\max].
\end{equation}

\subsection{Modeling the Minimum Impulse Bit}
\label{subsection:embedding:mib}

\begin{figure}
  \centering
  \begin{subfigure}[t]{0.33\textwidth}
    \includegraphics[page=1]{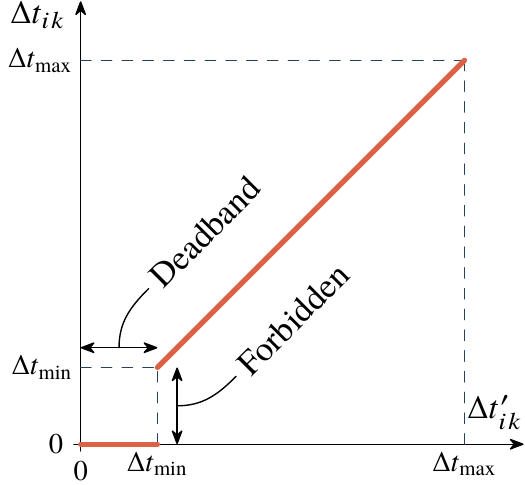}
    \caption{Exact deadband curve.}
    \label{fig:mib_exact}
  \end{subfigure}%
  \hfill%
  \begin{subfigure}[t]{0.33\textwidth}
    \includegraphics[page=2]{mib_smooth.pdf}
    \caption{Deadband approximation using \eqref{eq:smooth_discrete_logic}.}
    \label{fig:mib_homotopy_all}
  \end{subfigure}%
  \hfill%
  \begin{subfigure}[t]{0.33\textwidth}
    \includegraphics[page=3]{mib_smooth.pdf}
    \caption{Zoomed-in view of (\protect\subref{fig:mib_homotopy_all}) around
      $\Dt[\min]$.}
    \label{fig:mib_homotopy_zoom}
  \end{subfigure}
  \caption{Approximation of the MIB constraint \optieqref{brp}{mib} using
    multinomial logit smoothing \eqref{eq:smooth_discrete_logic}.}
  \label{fig:mib_illustration}
\end{figure}


The MIB constraint \optieqref{brp}{mib} is the most intricate one to model
effectively, and has been the core motivation behind developing a new way to
handle discrete logic constraints. Our past work used STCs, which exhibited
locking and prevented the algorithm from converging effectively in some cases
\cite{MalyutaScitechDocking}. Among the several possible ways of fitting the MIB
constraint into the discrete logic framework of \eqref{eq:discrete_logic_couple},
we present one way that yields good convergence performance across a wide
variety of instances of \pref{brp}.

Let us focus the discussion on pulse $\Dt[ik]$, in other words the \th{$i$}
thruster at the \th{$k$} control opportunity. We view the thruster as an
actuator with a deadband, as illustrated in \figref{mib_exact}. The ``input'' or
``reference'' pulse duration is given by a \textit{continuous} variable
$\Dtr[ik]\in[0,\Dt[\max]]$. When this value falls below $\Dt[\min]$, the
``obtained'' or ``output'' pulse duration which the thruster actually executes
is zero. Thus, while $\Dtr[ik]$ is a continuous variable that can take any
value in $[0,\Dt[\max]]$, the obtained pulse duration $\Dt[ik]$ exhibits a jump
discontinuity at $\Dt[\min]$. Modeling this jump discontinuity is precisely the
focus of our smooth approximation strategy.

Comparing \figref{mib_exact} with the standard model \eqref{eq:discrete_logic_couple},
we can write the following \textifelse logic:
\begin{subequations}
  \label{eq:mib_if_else_basic}
  \begin{align}
    \Dtr[ik]\le\Dt[\min]~&\implies~\Dt[ik]=0, \\
    \Dtr[ik]>\Dt[\min]~&\implies~\Dt[ik]=\Dtr[ik].
  \end{align}
\end{subequations}

We can thus define $\ng=1$, $z=\brak[big]{\Dt[ik]; \Dtr[ik]}$, and use the
predicate:
\begin{equation}
  \label{eq:mib_predicate}
  \predicate[1](\Dtr[ik]) = \Dtr[ik]-\Dt[\min].
\end{equation}

This predicate is used in \eqref{eq:discrete_logic_or_approx} to form
$\smoothRmib$, the smooth \textor approximation for the MIB predicate. As for
the implications on the right\dash hand side of \eqref{eq:mib_if_else_basic}, we
can use pairs of inequalities to represent equality constraints as required by
\eqref{eq:discrete_logic_couple}. This yields the following \textif and \textelse
implications:
\begin{equation}
  \label{eq:mib_if_else_implications}
  \fL(\Dt[ik]) = \Matrix{\Dt[ik] \\ -\Dt[ik]},\quad
  \fR(\Dt[ik], \Dtr[ik]) = \Matrix{\Dt[ik]-\Dtr[ik] \\ \Dtr[ik]-\Dt[ik]}.
\end{equation}

Just like for the approach cone and plume impingement constraints,
\eqref{eq:mib_if_else_implications} can now be substituted into
\eqref{eq:smooth_discrete_logic} to obtain a smooth approximation of the deadband
behavior in \figref{mib_exact}. Simplifying the result, we obtain the following
constraint:
\begin{equation}
  \label{eq:smooth_mib}
  \Dt[ik] = \smoothRmib(\Dtr[ik])\Dtr[ik].
\end{equation}

The smooth approximation is shown in \figref{mib_homotopy_all} for a number of
homotopy parameter $\sharp$ values. We call this approximation the smooth
deadband curve (SDC). As $\sharp$ increases, the approximation converges to the
exact deadband curve \textit{with one significant exception}: the ``forbidden''
region (i.e., the jump discontinuity) from \figref{mib_exact} becomes part of the
SDC as a quasi\dash vertical ``wall'' for large $\sharp$ in
\figref{mib_homotopy_all}. This raises the followig question: can a rendezvous
trajectory exploit this wall and therefore ``get around'' the MIB constraint?
Alas, the answer is yes, and our numerical tests show that this happens quite
regularly. Generally, this adversarial exploitation of the model feeds into a
long\dash standing pain point of optimization. As Betts writes in
\cite[p.~701]{BettsBook}, ``If there is a flaw in the problem formulation, the
optimization algorithm will find it.'' To fix this side effect and forbid
$\Dt[ik]$ from exploiting the wall, we introduce a new constraint to the
optimization problem. 

\subsubsection{The Wall Avoidance Constraint}
\label{subsubsection:embedding:mib:dont}

We now develop an extra constraint to ensure that no $\Dt[ik]$ can exploit the
wall part of the SDC \eqref{eq:smooth_mib}. We ask ourselves the following
question: what makes the wall different from the other parts of the SDC?  One
property stands out above all others: for large $\sharp$ values the wall has a
very large gradient, as opposed to other parts of the curve where the gradient
is approximately zero or one. There is another favorable property of
\eqref{eq:smooth_mib}: in the limit as $\sharp$ increases, the smooth approximation
converges to a function whose gradient monotonically increases for
$\Dt[ik]\in[0,\Dt[\min])$, and monotonically decreases for
$\Dt[ik]\in(\Dt[\min],\Dt[\max]]$. In other words, \eqref{eq:smooth_mib} has an
\textit{inflection point} at $\Dt[\min]$ for large $\sharp$, where its gradient
takes its maximum value. We call this the ``pivot'' since the SDC appears to
revolve around this point as $\sharp$ increases. This is visible in
\figref{mib_homotopy_all,mib_homotopy_zoom} for the brighter colored curves that
correspond to larger $\sharp$ values.

\defvar[hide]{$\Dtbuff$}{buffer zone around $\Dt[\min]$ for the wall avoidance
  constraint, \si{\second}}

We develop the following intuition from the above discussion: if we constrain
$\Dt[ik]$ such that the SDC's gradient is sufficiently less than its value at
the pivot, then $\Dt[ik]$ cannot exploit the wall. To put this into practice,
define $\Dtbuff$ to be a ``buffer'' around $\Dt[\min]$. We want the gradient at
$\Dt[ik]$ to be less than its value at the buffered pulse duration
$\Dt[\min]+\Dtbuff$. The SDC gradient at $\Dt[\min]+\Dtbuff$ is computed as
follows using \eqref{eq:smooth_mib}:
\begin{equation}
  \label{eq:deadband_gradient_buffer}
  \Gbuff \definedas
  \left.\frac{\dd\smoothRmib(\Dtr[ik])}{\dd\Dtr[ik]}
  \right|_{\Dtr[ik]=\Dt[\min]+\Dtbuff}
  (\Dt[\min]+\Dtbuff)+
  \smoothRmib(\Dt[\min]+\Dtbuff).
\end{equation}

This allows us to impose the following \textit{wall avoidance constraint},
which prevents $\Dt[ik]$ from taking values along the wall of the SDC:
\begin{equation}
  \label{eq:dont_do_that_constraint}
  \frac{\dd\smoothRmib(\Dtr[ik])}{\dd\Dtr[ik]}\Dtr[ik]+
  \smooth{\Rif}(\Dtr[ik]) \le \Gbuff.
\end{equation}

\figref{mib_homotopy_zoom} illustrates an example region of $\Dtr[ik]$ and
$\Dt[ik]$ values that is effectively removed by
\eqref{eq:dont_do_that_constraint}. In the figure, $\Dt[\min]=0.2~\si{\second}$ and
$\Dtbuff=0.06~\si{\second}$. The gradients of all points inside the red region
are larger than $\Gbuff$, hence the corresponding choices of $\Dt[ik]$ are
infeasible. Because the aforementioned monotonicity property guarantees that
this region contains the wall, the net effect is that the SDC wall can no
longer be exploited by the optimization.

\subsubsection{Improving Convergence}
\label{subsubsection:embedding:mib:convergence}

The smoothed MIB constraint \eqref{eq:smooth_mib} introduced a new input variable
$\Dtr[ik]$ to represent a reference pulse duration. This variable was necessary
to model the deadband curve in \figref{mib_exact}. If we compare the deadband curve
to the original MIB constraint \optieqref{brp}{mib}, we realize that the only
``useful'' parts of the curve in \figref{mib_exact} that we actually need are the
origin (i.e., $\brak{\Dt[ik]; \Dtr[ik]}=0$) and the continuous trace
$\Dt[ik]=\Dtr[ik]$ where $\Dt[ik]>\Dt[\min]$. In both cases, we have the simple
relationship $\Dt[ik]=\Dtr[ik]$. Our numerical experience shows that
encouraging this equality significantly improves the convergence process of the
algorithm in \sref{scp}. We do this by adding the following regularization term to
the original cost \optiobjref{brp}:
\begin{equation}
  \label{eq:equality_cost}
  \Jeq = \weq\Dt[\min]\inv\sum_{i=1}^{\nrcs}\sum_{k=1}^{\Nc}
  \norm[1]{\Dt[ik]-\Dtr[ik]},
\end{equation}
where $\weq>0$ is some small weight for the cost. We view \eqref{eq:equality_cost}
as penalizing the choice $\Dt[ik]\ne\Dtr[ik]$. The use of the one-norm
encourages sparsity in the number of $\Dt[ik]$ that violate the equality. This
choice traces back to theory from lasso regression, sparse signal recovery, and
basis pursuit to compute sparse solutions via one-norm regularization
\cite{BoydConvexBook}.

\subsection{Smoothed Rendezvous Problem}
\label{subsection:embedding:summary}

We are now in a position to restate \pref{brp} as a continuous optimization problem
by using the smoothed discrete logic constraints from the previous
sections. The process is straightforward: simply replace each discrete logic
constraint with its smooth approximation. We call the result the smooth
rendezvous problem (SRP), stated below.
\begin{optimization}[
  label={srp},
  variables={x, \Dt, \Dtr, \tf},
  objective={\Jfuel+\Jeq}]
  &\mrm[t]{Dynamics~\optieqref{brp}{position}-\optieqref{brp}{angular_velocity}}, \\
  \optilabel{pulse_bounds}
  &0\le \Dt[ik]\le \Dt[\max],~0\le \Dtr[ik]\le \Dt[\max], \\
  \optilabel{mib}
  &\Dt[ik] = \smoothRmib(\Dtr[ik])\Dtr[ik], \\
  \optilabel{mib_dont}
  &\frac{\dd\smoothRmib(\Dtr[ik])}{\dd\Dtr[ik]}\Dtr[ik]+
  \smooth{\Rif}(\Dtr[ik]) \le \Gbuff, \\
  \optilabel{plume}
  &\Dt[ik] \le \smoothRplume\pare[big]{p(\minusone{k}\tc)}\Dt[\max]
  ~\mrm[t]{for all}~i\in\Ifr, \\
  \optilabel{approach}
  &\cos(\angappch)-\pare[big]{1+\cos(\angappch)}\smoothRappch(p)-
  \ex[\frL]\T p\norm[2]{p}\inv\le 0, \\
  &\mrm[t]{Boundary conditions~\optieqref{brp}{bcs1}-\optieqref{brp}{bcs3}}.
\end{optimization}

The key difference between \pref{brp} and the new \pref{srp} is that the latter no
longer contains integer variables to solve. Instead, there is a single homotopy
parameter $\sharp$ that regulates how accurately the smoothed constraints
\optieqref{srp}{mib}, \optieqref{srp}{plume}, and \optieqref{srp}{approach} approximate the
original discrete logic. Thus, we have eliminated the third difficulty
mentioned in \ssref{formulation:summary} (i.e., the mixed\dash integer programming
aspect). However, we are now faced with solving a nonconvex optimization
problem, and there remains the question of how to set the value of $\sharp$. In
the next section we answer both questions using sequential convex programming
and numerical continuation.

\section{Sequential Convex Programming with Numerical Continuation}
\label{section:scp}

We now present a numerical optimization algorithm that solves \pref{srp}. This
algorithm combines two key methodologies: sequential convex programming (SCP)
and numerical continuationg. SCP is an iterative scheme designed to solve
\pref{srp} for a given value of $\sharp$. The \textit{raison d'{\^e}tre} for
numerical continuation is to greatly expand the region of convergence of
iterative schemes \cite{Watson1986}. Due to the vanishing gradient problem and
the very large gradients at the ``step'' transition points of discrete logic
(see, for example,
\figref{smooth_approximation_stages,smooth_comparison,mib_illustration}), SCP is
unlikely to converge if a large $\sharp$ value is used right away together with
an initial guess that is not already almost optimal \cite{CSM2021}. As a result,
numerical continuation is used to aid SCP convergence. This is done by
providing an algorithm to update $\sharp$ starting from a small value where the
smooth approximation is coarse, and increasing it until a large value where the
approximation attains the accuracy level requested by the user.

Our core contribution is to \textit{merge} these two methods. In other words,
the algorithm that we present is not SCP with a numerical continuation ``outer
loop''. Rather, the methods are run simultaneously, which is a novel feature of
the proposed algorithm. The numerical results in \sref{results} show that this can
dramatically decrease the total number of required iterations without
sacrificing optimality.

\begin{figure}
  \centering
  \includegraphics[width=0.9\textwidth]{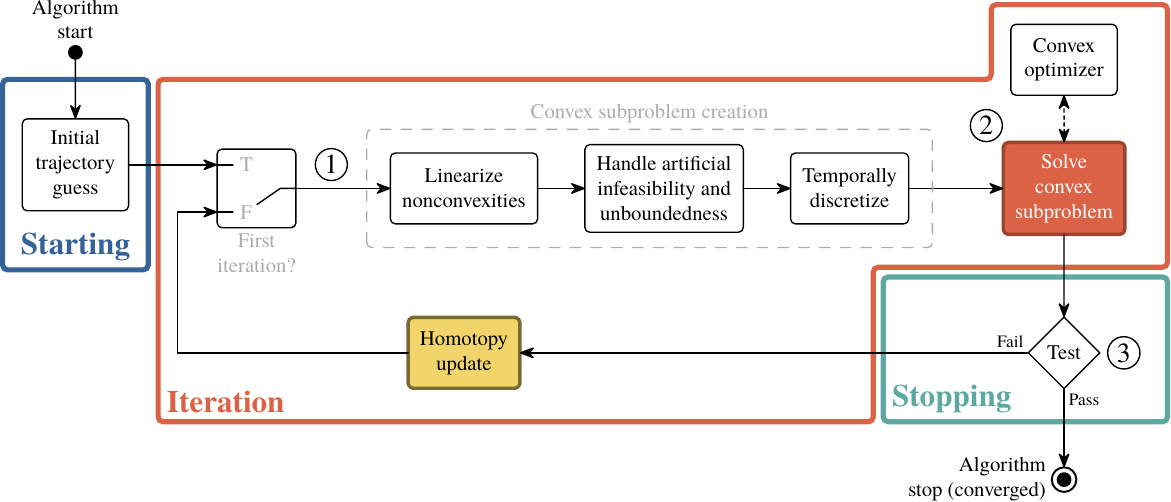}
  \caption{Block diagram illustration of the PTR sequential convex programming
    algorithm. Our novel contribution is the ``homotopy update'' block that
    implements numerical continuation to solve \pref{srp}.}
  \label{fig:scp_loop}
\end{figure}

\subsection{The Penalized Trust Region Algorithm}
\label{subsection:scp:ptr}

We begin by describing the penalized trust region (PTR) algorithm. This is a
particular SCP method that has been widely used for fast and even real\dash
time solution of nonconvex problems like \pref{srp}, where the value of
$\sharp$ is fixed \cite{Reynolds2020Real,SzmukThesis,ReynoldsThesis}. This
section provides a brief overview of PTR and identifies locations where the
method is changed in order to embed numerical continuation. These changes are
then described in the sections that follow. For the standard portions of the
PTR algorithm, we will refer the reader to existing literature which already
provides detailed explanations.

The goal of SCP in general, and PTR in particular, is to solve continuous\dash
time optimal control problems of the following form:
\begin{optimization}[%
  label={scp_general},%
  variables={x,u,p,\tf},%
  objective={J(x,u,p)}]%
  \optilabel{dynamics}%
  & \dot{x}(t) = f\pare[big]{t,x(t),u(t),p}, \\
  \optilabel{convex_path_constraints}%
  & \pare[big]{x(t),p} \in \set{X}(t),~\pare[big]{u(t),p} \in \set{U}(t), \\
  \optilabel{nonconvex_constraints}%
  & s\pare[big]{t,x(t),u(t),p} \leq 0, \\
  \optilabel{boundary_conditions}%
  & \gic\pare[big]{x(0),p} = 0,~\gtc\pare[big]{x(\tf),p} = 0,
\end{optimization}
where $x(\cdot)\in\reals^{\dimx}$ is the state trajectory,
$u(\cdot)\in\reals^{\dimu}$ is the control trajectory, and $p\in\real^{\dimp}$
is a vector of parameters. The function
$f : \real\times\real^{\dimx} \times \real^{\dimu} \times \real^{\dimp} \rightarrow \real^{\dimx}$
encodes the nonlinear equations of motion, which are assumed to be at least
once continuously differentiable. Initial and terminal boundary conditions are
enforced by using the continuously differentiable functions
$\gic:\real^{\dimx} \times \real^{\dimp} \rightarrow \real^{\dimgic}$ and
$\gtc : \real^{\dimx} \times \real^{\dimp} \rightarrow \real^{\dimgtc}$. The
convex and nonconvex path (i.e., state and control) constraints are imposed
using the convex sets $\set{X}(t)$, $\set{U}(t)$, and the continuously
differentiable function
$s : \real \times \real^{\dimx} \times \real^{\dimu} \times \real^{\dimp} \rightarrow \real^{\dimss}$. Finally,
a continuously differentiable cost function
$J:\reals^{\dimx}\times\reals^{\dimu}\times\reals^{\dimp}\to\reals$ encodes
some trajectory property that is to be minimized. Without giving the explicit
details here, we note that \pref{srp} can be fit into the mold of \pref{scp_general}
for any fixed value of $\sharp$. The interested reader may consult our
open\dash source implementation for details \cite{OurOpenSourceCode}, and may
refer to \cite{CSM2021} for a broad tutorial on the modeling process.

At the core of PTR is the idea of solving \pref{scp_general} through iterative
convex approximation. The algorithm can be represented in block diagram form as
shown in \figref{scp_loop}. The method is composed of three major parts: a way to
guess the initial trajectory ({\color{Blue}Starting}), an iteration scheme that
refines the trajectory until it is feasible and locally optimal
({\color{Red}Iteration}), and an exit criterion to stop once a trajectory has
been computed ({\color{Green}Stopping}). Strictly speaking, PTR is a nonlinear
local optimization algorithm known as a trust region method
\cite{ConnTrustRegionBook,NocedalBook,Kochenderfer2019}.

Let us begin by assuming that the homotopy parameter $\sharp$ is fixed to a
specific value. In other words, the ``homotopy update'' block in \figref{scp_loop}
is a simple feed\dash through that does nothing. PTR solves \pref{srp} using a
sequence of convex approximations called \textit{subproblems}. Roughly
speaking, the convex approximation is improved each time that a new solution is
obtained. Going around the loop of \figref{scp_loop}, all algorithms start with a
user-supplied initial guess, which can be very coarse (more on this later). At
\alglocation{\iterstartloc}, the SCP algorithm has available a so-called
reference trajectory, which may be infeasible with respect to the problem
dynamics and constraints. The nonconvexities of the problem are removed by a
local linearization around the reference trajectory, while convex parts of the
problem are kept unchanged. To ensure that linearization does not cause the
subproblems to become infeasible, extra terms are added which are known as
virtual controls (for the dynamics \optieqref{scp_general}{dynamics}) and virtual
buffers (for the constraints \optieqref{scp_general}{nonconvex_constraints} and
\optieqref{scp_general}{boundary_conditions}). The resulting convex continuous-time
subproblem is temporally discretized to yield a finite-dimensional convex
optimization problem. The optimal solution to the discretized subproblem is
computed at \alglocation{\solveloc}, where the SCP algorithm makes a call to
any appropriate convex optimization solver. The solution is tested at
\alglocation{\testloc} against stopping criteria. If the test passes, the
algorithm has converged and the most recent solution from
\alglocation{\solveloc} is returned. Otherwise, the solution becomes the new
reference trajectory for the next iteration of the algorithm.

The traditional PTR method as described above is covered in great depth in
existing literature. We refer the reader to a recent expansive tutorial
\cite{CSM2021}, and to papers which describe PTR in the context of rocket
landing, rendezvous and docking, and quadrotor flight
\cite{Szmuk2020,Reynolds2020,Reynolds2020Real,MalyutaScitechDocking,
  SzmukThesis,ReynoldsThesis}. In this paper we will focus our attention on the
novel ``homotopy update'' block in \figref{scp_loop}. This block implements a
numerical continuation method in order to update $\sharp$ until the smooth
approximations of discrete logic from \sref{embedding} become quasi\dash exact (in
other words, accurate to within a user-defined tolerance that can be
arbitrarily small).

\subsection{Non\dash embedded Numerical Continuation}
\label{subsection:scp:standard_continuation}

\begin{figure}
  \centering
  \begin{subfigure}[t]{0.49\textwidth}
    \centering
    \includegraphics[page=1,width=0.9\textwidth]{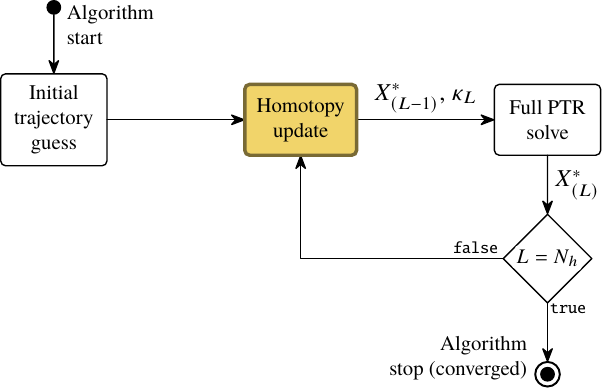}
    \caption{The non\dash embedded approach, where each homotopy update is
      followed by a full PTR solve.}
    \label{fig:numerical_continuation_standard}
  \end{subfigure}%
  \hfill%
  \begin{subfigure}[t]{0.49\textwidth}
    \centering
    \includegraphics[page=2,width=0.9\textwidth]{continuation_block_diagram.pdf}
    \caption{The proposed embedded approach, where homotopy updates occur
      between individual PTR iterations.}
    \label{fig:numerical_continuation_embedded}
  \end{subfigure}

  \caption{Comparison of the standard and embedded numerical continuation
    schemes. The ``Test'' in (\protect\subref{fig:numerical_continuation_embedded})
    corresponds to the stopping criterion from \figref{scp_loop}.}
  \label{fig:numerical_continuation_block_diagrams}
\end{figure}

In order to arrive at the embedded numerical continuation approach, we begin by
motivating a non\dash embedded scheme which we will then generalize to the
embedded algorithm. As shown in \figref{numerical_continuation_standard}, the basic
idea is to update the homotopy parameter $\sharp$ after each time that \pref{srp}
is solved for the current value of $\sharp$. Furthermore, each new call to PTR
is ``warm started'' by providing the most recent solution as the initial guess
in \figref{scp_loop}.

In formal terms, let $\Li$ denote the iteration number of the non\dash embedded
algorithm. Effectively, $\Li$ corresponds to the number of full PTR solves of
\pref{srp} that have occured up until the end of that iteration. If we place
ourselves at iteration $\Li$, then let $\sharp[\Li]$ denote the homotopy
parameter chosen by the ``homotopy update'' block, and let $\solution[\Li]$ be
the corresponding solution of \pref{srp} computed by PTR. Importantly, PTR is warm
started with the initial guess $\solution[\Li-1]$. When $\Li>1$, this
corresponds to the PTR solution from the previous iteration (i.e., the solution
of \pref{srp} for the previous value of $\sharp$). For the first iteration
$\Li=1$, $\solution[0]$ corresponds to the user\dash chosen initial trajectory
guess. The job of the homotopy update is the following: compute $\sharp[\Li]$
given $\solution[\Li-1]$ and $\sharp[\Li-1]$. While we describe the details
below, the basic idea is as follows: $\sharp[\Li]$ grows with $\Li$ and,
eventually, the smooth approximations from \sref{embedding} become quasi\dash exact
representations of the original discrete logic (e.g., see the example in
\figref{smooth_approximation_stages}). Once $\sharp[\Li]$ reaches some large
user\dash defined value that yields an accurate enough approximation of the
discrete logic, the algorithm terminates.

\begin{figure}
  \centering
  \begin{subfigure}[t]{0.49\textwidth}
    \centering
    \includegraphics[page=1]{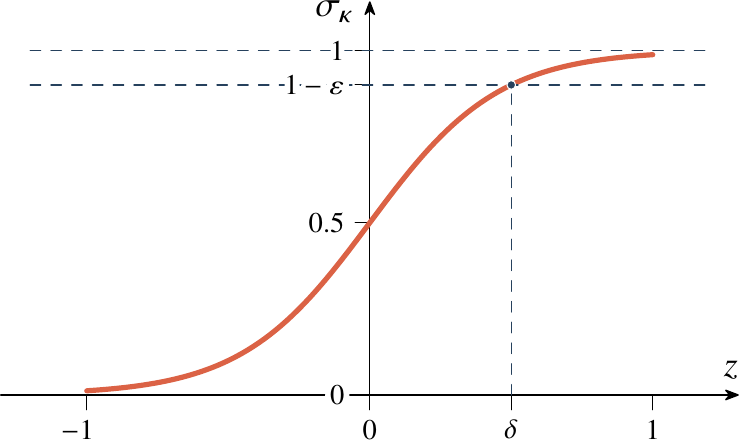}
    \caption{Each sigmoid is formed by fixing the values of $\sigeps$ and
      $\sigdelta$, which determine the homotopy parameter $\sharp$.}
    \label{fig:homotopy_update_annotations}
  \end{subfigure}%
  \hfill%
  \begin{subfigure}[t]{0.49\textwidth}
    \centering
    \includegraphics[page=2]{homotopy_update.pdf}
    \caption{The smoothness parameter $\sigdelta$ is decreased from
      $\sigdeltamax$ to $\sigdeltamin$. The sigmoid is sharpened as a result.}
    \label{fig:homotopy_update_sweep}
  \end{subfigure}
  \caption{The homotopy parameter $\sharp$ is updated by fixing a
    \textit{precision} $\sigeps$ and gradually reducing the \textit{smoothness}
    $\sigdelta$ where the sigmoid attains the value $1-\sigeps$.}
  \label{fig:homotopy_update}
\end{figure}

The remaining task for the non\dash embedded numerical continuation approach is
to define the internals of the homotopy update block in
\figref{numerical_continuation_standard}. Our method stems from viewing the sigmoid
function \eqref{eq:sigmoid} as a smooth model for a step function. As we increase
the homotopy parameter $\sharp$, we want to explicitly control how ``sharply''
the sigmoid approximates the step function's discontinuity. This leads us to
the following update rule, which is illustrated in \figref{homotopy_update}. As
shown in \figref{homotopy_update_annotations}, we define two parameters: a
\textit{precision} $\sigeps\in (0, 1)$ and a \textit{smootheness}
$\sigdelta>0$. The sigmoid function $\sigmoid[\sharp]$ in \eqref{eq:sigmoid} is
then required to satisfiy the following \textit{interpolation condition}: it
must equal $1-\sigeps$ when its argument equals $\sigdelta$. An exact step
function corresponds to $\sigeps=0$ and $\delta=0$, so we view $\sigeps$ and
$\sigdelta$ as defining how much the sigmoid deviates from the exact step
function.

For the homotopy update rule, we hold $\sigeps$ constant and define two bounds
on $\sigdelta$: a ``smoothest'' value $\sigdeltamax$ and a ``sharpest'' value
$\sigdeltamin<\sigdeltamax$. We then sweep $\sigdelta$ according to a geometric
progression:
\begin{equation}
  \label{eq:smootheness_sweep}
  \delta_{\siginterp} = \sigcoeff^{\siginterp}\sigdeltamax,~
  \quad~\sigcoeff = \sigdeltamin/\sigdeltamax,
\end{equation}
where $\siginterp\in [0,1]$ is an interpolation parameter. The effect is that
the sigmoid function is sharpened, as shown in \figref{homotopy_update_sweep}. The
homotopy value that satisfies the interpolation condition is given by:
\begin{equation}
  \label{eq:homotopy_value_interpolation}
  \sharp[\siginterp] =
  \frac{\ln\pare[big]{\sigeps\inv-1}}{\sigcoeff^{\siginterp}\sigdeltamax}.
\end{equation}

Equation \eqref{eq:homotopy_value_interpolation} defines a continuous range of
homotopy values from the smoothest ($\siginterp=0$) to the sharpest
($\siginterp=1$) case. In practice, we set a fixed number of updates $\Nhom$
and let $\alpha=(\Li-1)/(\Nhom-1)$ for $\Li=1,2,\dots,\Nhom$. Thus, $\Nhom$
defines the number of iterations in the non\dash embedded numerical
continuation algorithm of \figref{numerical_continuation_standard}. By substituting
this expression for $\alpha$ into \eqref{eq:homotopy_value_interpolation}, we
obtain the following formula for the homotopy value $\sharp[\Li]$ at iteration
$\Li$:
\begin{equation}
  \label{eq:homotopy_update_rule}
  \sharp[\Li] =
  \frac{\ln\pare[big]{\sigeps\inv-1}}{\sigcoeff^{(\Li-1)/(\Nhom-1)}\sigdeltamax}.
\end{equation}

\subsection{Embedded Numerical Continuation}
\label{subsection:scp:embedded_continuation}

We are now ready to describe the embedded numerical continuation algorithm
shown in \figref{numerical_continuation_embedded}. One key difference distinguishes
this algorithm from the non\dash embedded approach: PTR does not have to run to
completion before the homotopy parameter $\sharp$ is increased. As shown in
\figref{numerical_continuation_embedded}, the full PTR solve of the non\dash
embedded method is replaced by a \textit{single} PTR iteration (which
corresponds to the top half of the PTR block diagram in \figref{scp_loop}). We use
$\li$ to denote the PTR iteration counter. At each iteration $\li$, a homotopy
update rule is called that potentially changes the value of $\sharp$. This new
value and the most recent PTR iterate (i.e., subproblem solution) are used for
the next PTR iteration. The process shown in
\figref{numerical_continuation_embedded} works exactly like in \figref{scp_loop}, with
the blocks rearranged.

Now that we understand how the algorithm is structured, we need to describe the
homotopy update. This is composed of two parts: deciding \textit{whether} to
update $\sharp$, and then updating it. The latter piece works just like in the
previous section. Once we know that $\sharp$ should be updated, we use
\eqref{eq:homotopy_update_rule} to compute its new value:
\begin{equation}
  \label{eq:homotopy_update_rule}
  \sharp[\li] =
  \frac{\ln\pare[big]{\sigeps\inv-1}}{\sigcoeff^{\Li/(\Nhom-1)}\sigdeltamax},
  \quad \Li\gets \Li+1,
\end{equation}
where $\Li$ now represents the number of times that the homotopy parameter has
been updated so far (the count starts at $\Li=0$). The core of the embedded
homotopy update rule is the first piece: deciding whether to update
$\sharp$. For this, let $\cost[\li]$ denote the subproblem cost achieved at PTR
iteration $\li$. If the following condition holds, then we update $\sharp$:
\begin{equation}
  \label{eq:update_decision}
  \hombl\le\frac{\cost[\li-1]-\cost[\li]}{\abso{\cost[\li-1]}}\le\hombu\quad
  \land\quad L<\Nhom.
\end{equation}

The second half of the condition is simple: don't update $\sharp$ if this is
already its highest value. The first half is a condition on relative cost
decrease over the past iteration. If the cost in the current iteration
decreased by less than $\hombu$ relative to the last iteration, then the
algorithm is ``converging'' for the current value of $\sharp$ and it is time to
update it. However, the cost is not guaranteed to decrease monotonically with
PTR iterations. Thus, the relative cost decrease may be negative, which means
that the cost increased over the past iteration. In this case, we may specify a
certain (small) tolerance $\hombl<0$. This means that we will still update
$\sharp$ if the cost did not increase by more than $\hombl$ allows. In the
numerical results of \sref{results} we set $\hombl=-10^{-3}$ (i.e., a
$0.1\si{\percent}$ tolerance).

\begin{algorithm}
  \begin{algorithmic}[1]
    \State $X_{(0)}\gets\textnormal{initial trajectory guess}$,~%
    $L\gets 0$,~$\li\gets 0$
    \While{\texttrue}
    \State $\li\gets\li+1$
    \If{$\li=1$ or \eqref{eq:update_decision} is \texttrue}
    \label{alg:embedded_numerical_continuation:line:hom_start}
    \State Execute the update rule \eqref{eq:homotopy_update_rule}
    \Else
    \State $\sharp[\li]\gets\sharp[\li-1]$
    \EndIf
    \label{alg:embedded_numerical_continuation:line:hom_end}
    \State $X_{(\li)}\gets\textnormal{do one PTR step (see \figref{scp_loop}) using
      $\sharp[\li]$ and the reference trajectory $X_{(\li-1)}$}$
    \If{$\Li=\Nhom$ and $X_{(\li)}$ passes the Test at
      location \alglocation{\testloc} in \figref{scp_loop}}
    \label{alg:embedded_numerical_continuation:line:stopping}
    \State \textbf{return} $X_{(\ell)}$ \Comment{Converged}
    \EndIf
    \EndWhile
  \end{algorithmic}
  \caption{The proposed sequential convex programming algorithm with embedded
    numerical continuation. The method can solve optimal control problems with
    discrete logic constraints.}
  \label{alg:embedded_numerical_continuation}
\end{algorithm}

To summarize the above discussion, \algref{embedded_numerical_continuation}
formalizes the embedded numerical continuation scheme illustrated in
\figref{numerical_continuation_embedded}. On
\algref[start=hom_start,end=hom_end]{embedded_numerical_continuation}, which
correspond to the ``homotopy update'' block of
\figref{numerical_continuation_embedded}, a decision is made using
\eqref{eq:update_decision} whether to update the current $\sharp$ value. If the
answer is \texttrue, then $\sharp$ is updated using
\eqref{eq:homotopy_update_rule}. Otherwise, it is maintained at its present
value. The algorithm iterates in this way until the stopping criterion on
\algref[start=stopping]{embedded_numerical_continuation} is satisfied.

\section{Numerical Results}
\label{section:results}

\begin{figure}
  \centering
  \includegraphics{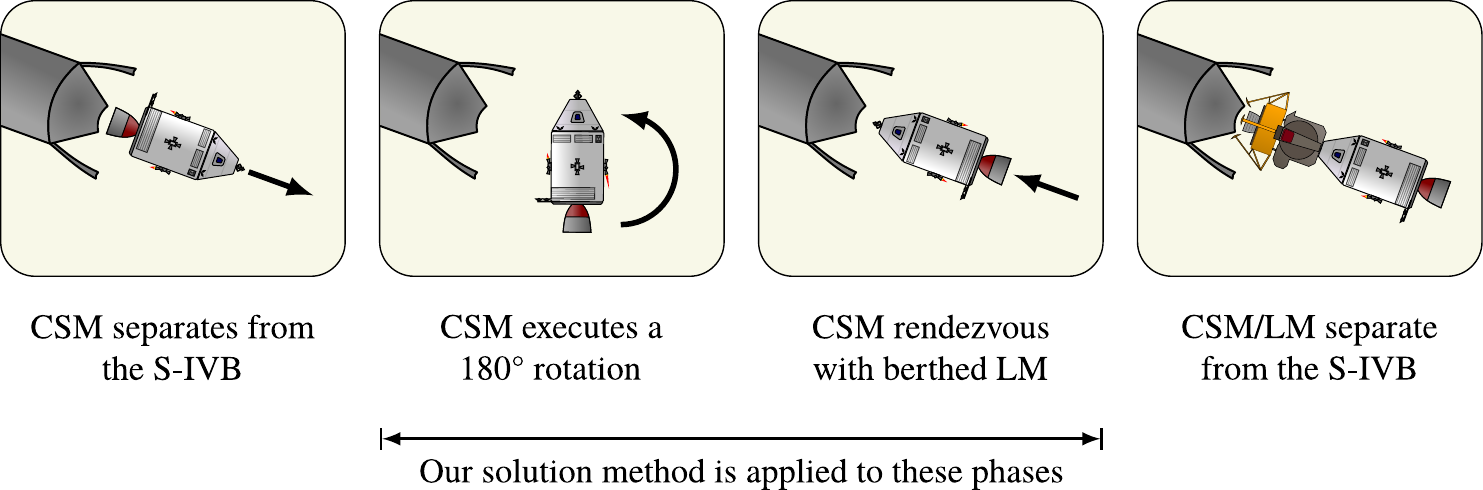}
  \caption{Illustration of the Apollo CSM Transposition and Docking maneuver
    with the LM housed inside the Saturn S-IVB third stage
    \cite[Figure~2-11]{NASA_apollo_mass}.}
  \label{fig:transposition_docking_illustration}
\end{figure}

In this section we apply our solution method to a more challenging variant of
the historical docking maneuver between the Apollo Command and Service Module
(CSM) and the Lunar Module (LM). \ssref{results:parameters} defines the problem
parameters and \ssref{results:plots} discusses the solved trajectory and various
computational aspects. The key takeaways are as follows. Our algorithm is able
to consistently find rendezvous trajectories that satisfy the discrete logic
constraints from \sref{formulation}. The algorithm is insensitive to the
$\hombu$ tolerance parameter in \eqref{eq:update_decision}, and in fact increasing
this value can dramatically reduce the total number of iterations. The total
convex solver time is approximately $13.5~\si{\second}$, which is fast for an
implementation that is not optimized for speed. \ssref{results:plots} discusses how
the algorithm can be made to run in under $10~\si{\second}$ of total solution
time.

\subsection{Problem Parameters}
\label{subsection:results:parameters}

The numerical example is inspired by the Apollo CSM ``Transposition and
Docking'' (TD) maneuver \cite[Section~2.13.1.1]{NASA_aoh_vol1}. As illustrated
in \figref{transposition_docking_illustration}, this maneuver uses the RCS
thrusters of the CSM in order to dock with the LM, which is housed inside the
Saturn S-IVB third stage. The most interesting feature of this maneuver is the
$180\si{\degree}$ rotation of the CSM, which allows us to stress\dash test our
solution method for the CSM's coupled translation and rotation dynamics from
\ssref{formulation:chaser_dynamics}.

\begin{figure}
  \centering
  \includegraphics[width=0.45\textwidth]{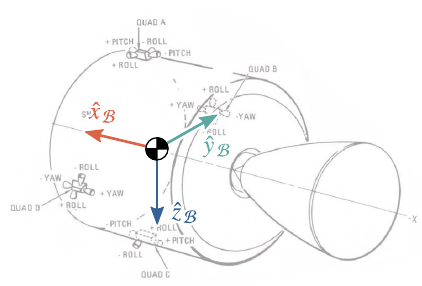}
  \caption{Layout of the Apollo SM RCS thrusters
    \cite[Figure~2.5-1]{NASA_aoh_vol1}. There are $\nrcs=16$ thrusters grouped
    into four ``quads'' labeled A/B/C/D and each having four independent
    hypergolic pressure-fed thrusters.}
  \label{fig:sm_rcs}
\end{figure}

The SM RCS system is composed of four similar, independent ``quads'' located
$90\si{\degree}$ apart around the Service Module (SM) circumference, as
illustrated in \figref{sm_rcs}. Each quad is composed of four independent
hypergolic pressure\dash fed pulse\dash modulated thrusters, yielding a total
of $\nrcs=16$ control inputs \cite[Section~2.5.1]{NASA_aoh_vol1}. We model the
complete, high\dash fidelity CSM geometry according to public NASA
documentation \cite{NASA_aoh_vol1,NASA_apollo_mass}. The CSM mass and inertia
are specified in \cite[Table~3.1-2]{NASA_apollo_mass}:
\begin{equation}
  m \approx 30323~\si{\kilo\gram},\quad
  J \approx
  \begin{bmatrix}
    49249  &   2862 &   -370 \\
    2862 & 108514  &  -3075 \\
    -370 &  -3075  & 110772
  \end{bmatrix}~\si{\kilo\gram\meter\squared}.
\end{equation}

The RCS thrusters are capable of producing approximately $\|\hat f_i\|_2=445$~N
of thrust in steady-state operation
\cite[Figure~2.5-8,~Table~4.3-1]{NASA_aoh_vol1}. However, during transposition
and docking they are pulse-fired in bursts
\cite[Section~2.5.1]{NASA_aoh_vol1,NASA_apollo_news}. In this mode of operation,
the minimum electric on-off pulse width is 12~ms
\cite[Section~2.5.2.3.1]{NASA_aoh_vol1}, generating an irregular burst of thrust
that lasts for upwards of 65~ms and with a peak of 300 to 350~N
\cite[Figure~2.5-9]{NASA_aoh_vol1}. To buffer the thrust away from this
irregular region, we set $\Delta t_{\min}=112~\si{\milli\second}$ (which
corresponds to a $50~\si{\newton\second}$ impulse) and
$\Delta t_{\max}=1~\si{\second}$. On a system architecture level, we assume
that irregularity in the thrust profile is going be corrected by a feedback
control system that tracks our open-loop rendezvous trajectory.

\newcolumntype{L}[1]{>{\raggedright\arraybackslash}p{#1}}
\newcolumntype{R}[1]{>{\raggedleft\arraybackslash}p{#1}}
\newcolumntype{C}[1]{>{\centering\arraybackslash}p{#1}}
\begin{table}[t]
  \centering

  \renewcommand*{\arraystretch}{1.2}

  \newcommand{\grouptitle}[1]{%
    \multicolumn{2}{c}{\itshape #1}}

  \begin{subfigure}[b]{0.48\textwidth}
    \begin{tabular}{C{2.6cm}L{4.3cm}}
      \hline\hline
      \multicolumn{1}{c}{Parameter}
      & \multicolumn{1}{l}{Value} \\
      \hline
      $\Dt[\min]$ & $112~\si{\milli\second}$ \\
      $\Dt[\max]$ & $1~\si{\second}$ \\
      $\Frcs$ & $445~\si{\newton}$ \\
      $\rplume$ & $20~\si{\meter}$ \\
      $\rappch$ & $30~\si{\meter}$ \\
      $\angappch$ & $10\si{\degree}$ \\
      $\tf$ & $\in [100,1000]~\si{\second}$ \\
      $\norb$ & $400~\si{\kilo\meter}$ LEO \\
      $p_0$ & $100\ex[\frL]-20\ez[\frL]+20\ey[\frL]~\si{\meter}$ \\
      $v_0$ & $0~\si{\meter\per\second}$ \\
      $q_0$ & $\brak{0;0;0;1}$ \\
      $\omega_0,~\omega_f$ & $0~\si{\radian\per\second}$ \\
      $p_f$ & $4.48\ex[\frL]-0.05\ey[\frL]+0.17\ez[\frL]~\si{\meter}$ \\
      $v_f$ & $-0.1\ex[\frL]~\si{\meter\per\second}$ \\
      $q_f$ & $\brak{0;0.26;0.97;0}$ \\
      $\tol{r_f}$ & $0.1~\si{\meter}$ \\
      $\tol{v_f}$ & $1~\si{\centi\meter\per\second}$ \\
      $\tol{q_f}$ & $1\si{\degree}$ \\
      $\tol{\omega_f}$ & $0.01\si{\degree\per\second}$ \\
      \hline \hline
    \end{tabular}
    \caption{Rendezvous parameters.}
  \end{subfigure}%
  \hfill
  \begin{subfigure}[b]{0.48\textwidth}
    \begin{tabular}{C{2.6cm}L{4.3cm}}
      \hline\hline
      \multicolumn{1}{c}{Parameter}
      & \multicolumn{1}{l}{Value} \\
      \hline
      $\sigeps$ & $10^{-2}$ \\
      $\sigdeltamax$ & $10$ \\
      $\sigdeltamin$ & $0.01$ \\
      $\Nhom$ & $10$ \\
      $\weq$ & $1$ \\
      $\hombl$ & $-10^{-3}$ \\
      $\hombu$ & $0.1$ \\
      $\Dtbuff$ & $11.2~\si{\milli\second}$ \\
      \hline \hline
    \end{tabular}
    \caption{Algorithm parameters.}
  \end{subfigure}%
  \vspace{1mm}

  \caption{Numerical parameters for the Apollo CSM/LM Transposition and
    Docking trajectory.}
  \label{table:parameters}
\end{table}

\tabref{parameters} summarizes the major numerical values used to obtain the
results of the following section. Other parameters not mentioned (such as the
CSM geometry) can be consulted directly in our open\dash source implementation
\cite{OurOpenSourceCode}. Note that the maneuver we are seeking to optimize is
more complicated than the original Apollo TD concept of operations. The Apollo
initial position $p_0$ was almost purely along the $\ex[\frL]$ axis, whereas we
add significant $\ey[\frL]$ and $\ez[\frL]$ displacement in order to stress the
algorithm. Furthermore, the original TD maneuver takes place after translunar
injection whereas we assume a circular lower Earth orbit. This allows us to use
the Clohessy-Wiltshire-Hill dynamics \eqref{eq:eom_lvlh_relative_motion}, which
adds further complexity compared to our previous work
\cite{MalyutaScitechDocking}.

Our algorithm from \sref{scp} is implemented using the framework introduced in
\cite{CSM2021}. The Julia programming language is used because it is simple to
read like Python, yet it can be as fast as C/C++ \cite{bezanson2017julia}. The
timing results in the next section correspond to a Dell XPS 13 9260 laptop
powered by an Intel Core i5-7200U CPU clocked at 2.5~\si{\giga\hertz}. The
computer has 8~\si{\gibi\byte} LPDD3 RAM and 128~\si{\kibi\byte} L1,
512~\si{\kibi\byte} L2, and 3~\si{\mebi\byte} L3 cache. The ECOS software
(written in C) is used as the low\dash level numerical convex solver at
location \alglocation{\solveloc} in \figref{scp_loop} \cite{domahidi2013ecos}.

\subsection{Computed Trajectory}
\label{subsection:results:plots}

\begin{figure}
  \centering
  \includegraphics[width=\textwidth]{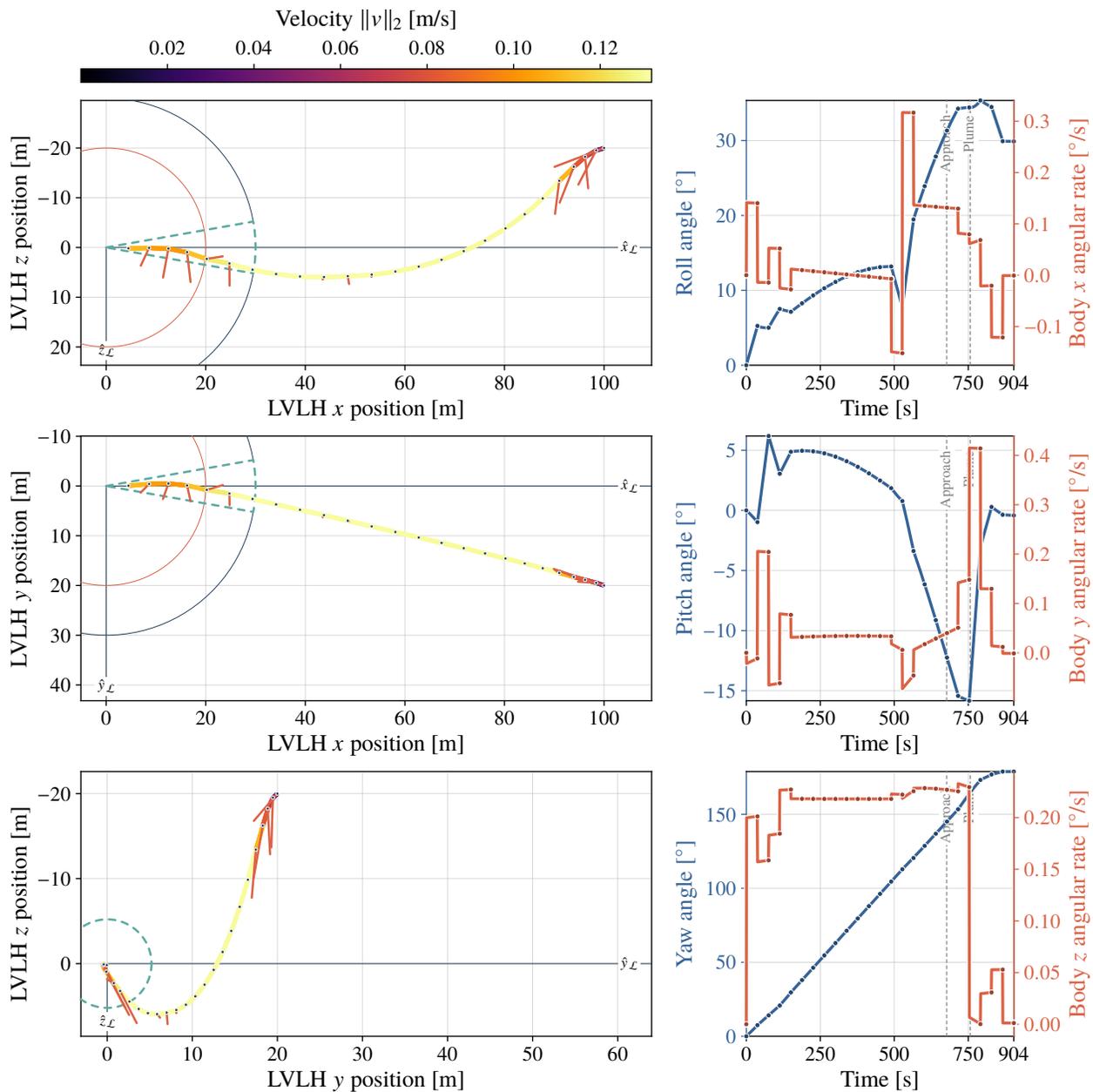}
  \caption{Projected views of the optimized trajectory in the LVLH frame (left)
    and the corresponding evolution of the chaser's attitude (right).%
  }
  \label{fig:trajectory_2d}
\end{figure}

\begin{figure}
  \centering
  \includegraphics[width=\textwidth]{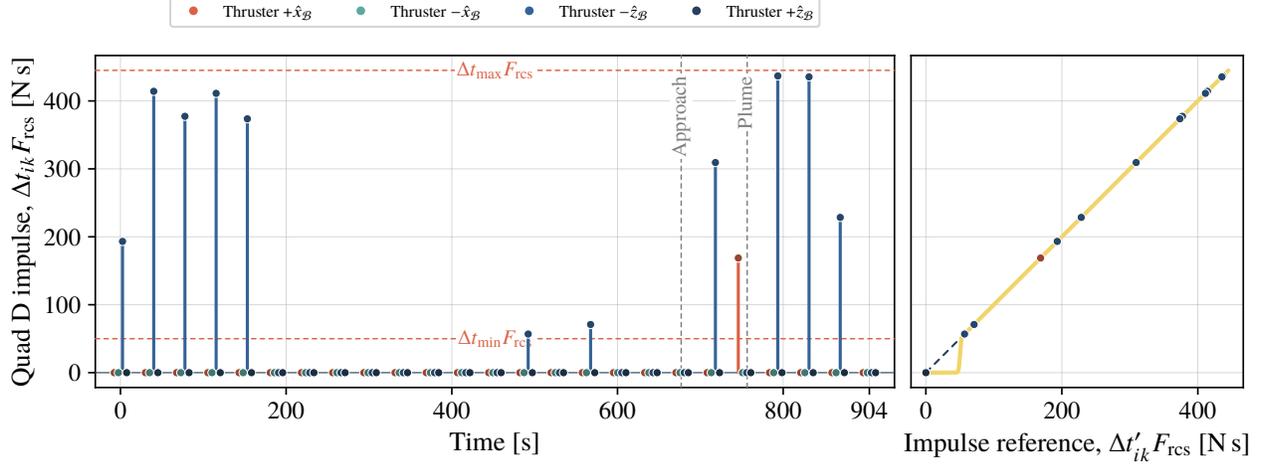}
  \caption{RCS thruster firing history for quad D of \figref{sm_rcs}, shown in
    terms of the impulse $\Frcs\Dt[ik]$. In the left plot, the four thrusters
    are ``spread'' along the time axis in order to not overlap.}
  \label{fig:inputs}
\end{figure}

\begin{figure}[t]
  \centering
  \includegraphics[width=1\textwidth]{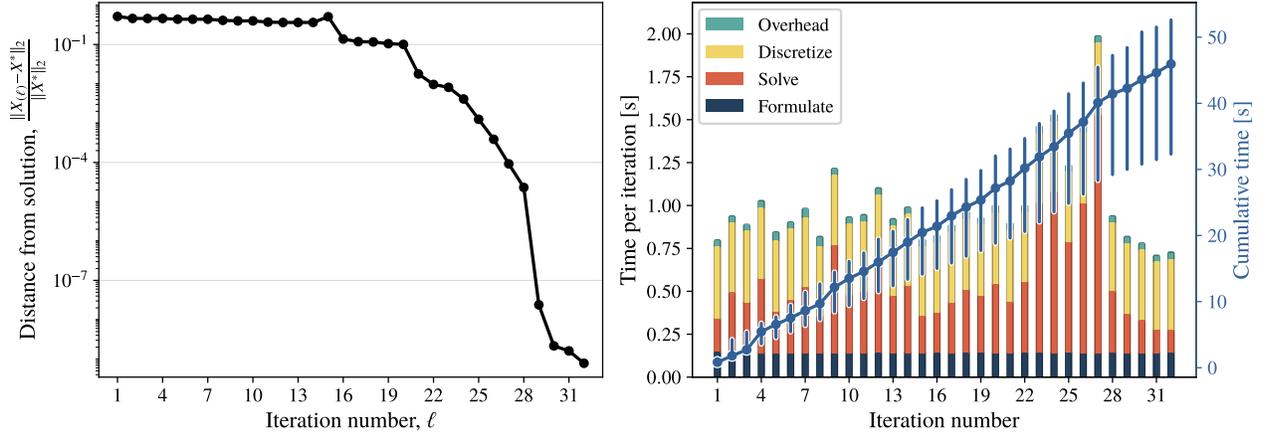}
  \caption{Convergence and runtime performance of our algorithm. The runtimes
    show statistics over 20 executions. Median values are shown, and error bars
    show the $10\%$ (bottom) and $90\%$ (top) quantiles.}
  \label{fig:convergence}
\end{figure}

\begin{figure}[t]
  \centering
  \includegraphics[width=\textwidth]{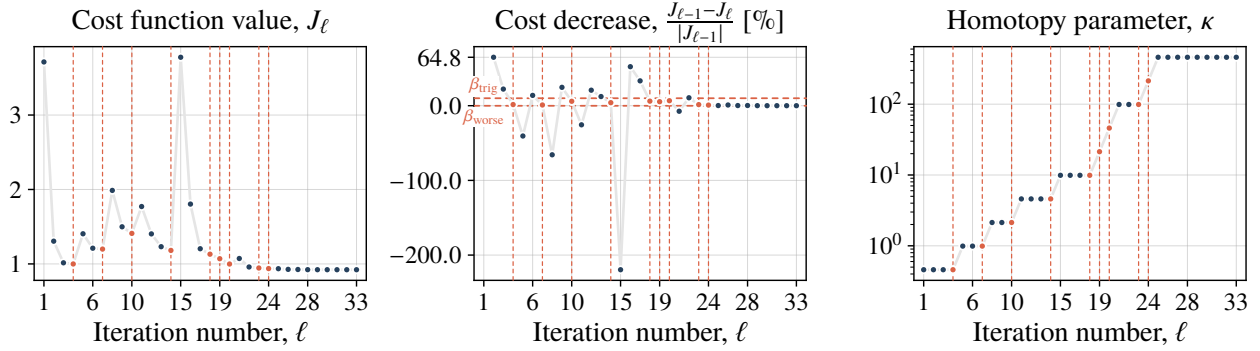}
  \caption{%
    Evaluation of the cost function \optiobjref{srp} over the PTR
    iterations (note the log scale for $\sharp$). Homotopy parameter updates
    are highlighted by {\color{Red}red} dashed vertical
    lines and markers.
  }
  \label{fig:cost}
\end{figure}

\begin{figure}
  \centering
  \includegraphics[width=0.8\textwidth]{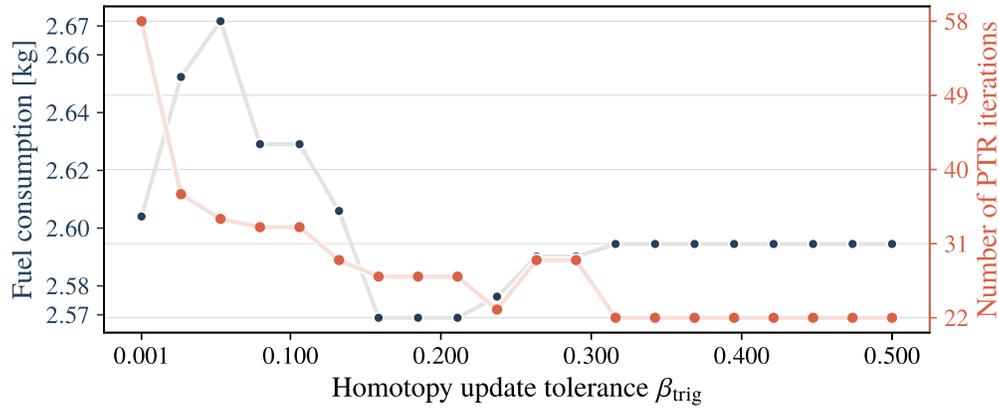}
  \caption{Dependence of the converged trajectory's fuel consumption and of our
    algorithm's total iteration count on the value of $\hombu$ in
    \eqref{eq:update_decision}.}
  \label{fig:homotopy_threshold}
\end{figure}

Figures~\ref{fig:trajectory_2d}-\ref{fig:homotopy_threshold} exhibit our algorithm's
solution as well as its associated properties for \pref{srp} with the parameters in
\tabref{parameters}. The initial guess provided to the algorithm in \figref{scp_loop}
is a straight\dash line interpolation in position and a spherical linear
interpolation for the attitude quaternion \cite{Sola2017}. The initial RCS
thruster pulse durations are all set to zero.

We begin by discussing the position trajectory, which is shown in the LVLH
frame in the left column of \figref{trajectory_2d}. Projections are also shown for
the approach sphere ({\color{DarkBlue}blue}), the approach cone
({\color{Green}green}), and the plume impingement sphere
({\color{Red}red}). The {\color{Red}red} vectors represent the direction and
relative magnitude of the net thrust generated by the combined action of the
RCS thrusters. The circular markers show the chaser's COM for the discrete\dash
time solution, while the continuous trajectory is obtained by integrating the
optimal control through the original nonlinear dynamics of
\ssref{formulation:chaser_dynamics}. Because the two trajectories coincide, we
conclude that the converged trajectory is dynamically feasible.

The trajectory in \figref{trajectory_2d} has two salient features. First, the RCS
thrusters fire mostly at the start to initiate motion, and near the end to
cancel vehicle rates just prior to docking. This resembles the classical
two\dash impulse rendezvous maneuver \cite{CurtisBook}, modified to account for
6-DOF dynamics, the RCS system geometry, and the discrete logic constraints
\optieqref{srp}{pulse_bounds}-\optieqref{srp}{approach}, all of which are absent in the
classical setup. Secondly, recall that negative $\ez[\frL]$ positions
correspond to lower orbits where objects move faster relative to the
target. The chaser exploits this ``gift'' from orbital mechanics by dipping
into the negative $\ez[\frL]$ positions (see the top and bottom plots) where it
benefits from a zero-fuel acceleration to the target. Furthermore, note how the
chaser stays within the approach cone when it is inside the approach sphere, as
required by \optieqref{srp}{approach}.

The evolution of the chaser's attitude along this trajectory is shown in the
right column of \figref{trajectory_2d}. The quaternion attitude was converted to
the more informative Euler angles using the Tait-Bryan yaw\dash pitch\dash roll
convention. {\color{Green}Green} vertical lines demarcate the times at which
the chaser enters the approach and plume imbpingement spheres. Velocity and
angular rate states exhibit jumps according to our impulsive thruster model in
\ssref{formulation:impulsive_thrust}. Note that the chaser assumes a
$30\si{\degree}$ roll angle at docking, as required by the CSM/LM geometry
\cite[Figure~2-4]{NASA_apollo_mass}.

The RCS thruster pulse history is shown in \figref{inputs} for quad D from
\figref{sm_rcs}, which is representative of the pulse histories for the other
quads. The pulses are relatively sparse and clustered around the start and end
of the trajectory. As required by the plume impingement constraint
\optieqref{srp}{plume}, the forward thrusters are silent once the chaser is inside
the plume impingement sphere. Furthermore, some pulse durations are almost
exactly $\approx\Dt[\min]~\si{\second}$. This shows that the smoothed discrete
logic \optieqref{srp}{mib} actively enforces the MIB constraint
\eqref{eq:pulse_duration_constraint}. The constraint \optieqref{srp}{mib} is
indispensable for satisfying the minimum impulse\dash bit, and removing it
causes the MIB constraint to be violated.

We can estimate the total fuel consumption of the rendezvous trajectory using
NASA charts for RCS thruster performance
\cite[Figure~4.3-6~and~4.3-7]{NASA_csm_aoh}. These charts map pulse duration to
the corresponding amount of fuel consumed by a single thruster. By applying
these charts to the pulse history in \figref{inputs}, we obtain a fuel consumption
of $2.63~\si{\kilo\gram}$. Unfortunately, NASA documentation on the actual fuel
consumption achieved by the Apollo missions is unclear;
\cite[Table~3.1-7]{NASA_apollo_mass} suggests that it was $32~\si{\kilo\gram}$,
but this confounds the other phases of the TD maneuver which we do not consider
(see \figref{transposition_docking_illustration}). In any case, it appears that our
trajectory uses considerably less fuel, not to mention that its initial
conditions are more challenging than those of the Apollo concept of operations
due to the initial position offsets along $\ey[\frL]$ and $\ez[\frL]$.

The convergence process of our algorithm and the runtime performance of its
implementation are shown in \figref{convergence}. {\color{DarkBlue}Formulate}
measures the time taken to parse the subproblem into the input format of the
convex optimizer; {\color{Red}Solve} measures the time taken by the core convex
numerical optimizer; {\color{Yellow}Discretize} measures the time taken to
temporally discretize the linearized dynamics from
\ssref{formulation:chaser_dynamics}; and {\color{Green}Overhead} measures the time
taken by all other supporting tasks during a single PTR iteration.

The algorithm appears to attain a superlinear convergence rate (noticeable over
iterations $\ell\in [19,30]$). A small spike in solver time appears around the
iterations where the homotopy parameter changes rapidly (see \figref{cost}
ahead). Otherwise, the subproblem difficulty stays roughly constant over the
iterations. While our Julia implementation takes a median time of
$50~\si{\second}$, the cumulative median time for solving the subproblems at
location \alglocation{\solveloc} in \figref{scp_loop} is $13.5~\si{\second}$ (which
is the sum of the {\color{Red}Solve} bars in \figref{convergence}). This
corresponds to the time taken by the ECOS convex solver, which is written in
C. We view this as a favorable runtime result for the following reasons, which
we state based on experience from \cite{Reynolds2020Real}. ECOS is a generic
solver, and a custom solver is likely to run at least twice as fast
\cite{Dueri2014,Dueri2017}. Coupled with other implementation efficiencies, we
expect that the total solver time can be reduced to
$<5~\si{\second}$. Furthermore, our code is optimized for readability. By
writing other parts of the algorithm in a compiled language and optimizing for
speed, we can expect to shrink the other $36~\si{\second}$ of runtime down to
$<5~\si{\second}$ as well. Thus, a speed-optimized implementation of our
algorithm can likely solve the rendezvous problem in under $10~\si{\second}$,
which is quite acceptable for rendezvous applications since the actual
trajectory can last for several thousand seconds.

\figref{cost} shows the evolution of the cost function value over the PTR
iterations. Every time the cost improvement falls within the decision range of
\eqref{eq:update_decision}, the homotopy parameter is updated. The update is
followed by a spike in the cost, with fast subsequent improvement to an equal
or better (i.e., smaller) value. During the final stages of the optimization
(iterations $\ell\ge 18$), increases in $\sharp$ no longer cause appreciable
spikes in cost. This is remarkable, given that it is over these iterations that
the homotopy parameter experiences its largest growth (since it grows
exponentially, as seen in \figref{homotopy_update_sweep} and the log scale of the
rightmost plot in \figref{cost}). This means that, well before convergence occurs,
our algorithm already finds a solution that is feasible with respect to the
final ``sharp'' approximation of the discrete logic. This analysis is
corroborated by the left plot in \figref{convergence}, where it can be seen that
past iteration $\ell\approx 20$ the amount by which the solution changes drops
off quickly.

Finally, \figref{homotopy_threshold} analyzes the depence of the optimal solution
and of our algorithm's performance on the choice of homotopy update tolerance
$\hombu$ in \eqref{eq:update_decision}. This reveals several favorable properties
of the algorithm. First, by increasing $\hombu$ we can dramatically lower the
total iteration count and speed up the solution time. A very low value of
$\hombu$ emulates the non\dash embedded numerical continuation scheme from
\figref{numerical_continuation_standard}, since $\sharp$ does not update until PTR
has quasi\dash converged for its current value. By increasing $\hombu$, we can
lower the iteration count by over $60\si{\percent}$ for this rendezvous
example. We observe this behavior consistenly across different initial
conditions. At the same time as lowering the iteration count, we basically
maintain a consistent level of fuel\dash optimality. The fuel consumption goes
up and down slightly, but on balance there is no perceptible trend. A notable
downside of using a larger $\hombu$ is an increased danger of not converging to
a feasible trajectory, since we have ``rushed'' the algorithm too much. This
does not happen in the present example, but we have noticed the issue for other
initial conditions. Our future work plans to investigate what is the
theoretically safe upper bound for the $\hombu$ value.

\section{Conclusion}
\label{section:conclusion}

This paper presents a novel algorithm combining sequential convex programming
with numerical continuation to handle a general class of discrete logic
constraints in a continuous optimization framework. This makes the approach
amenable to fast and reliable solution methods for trajectory optimization
problems commonly encountered in spaceflight. The algorithm is applied to the
terminal phase of rendezvous and docking maneuver, where a chaser spacecraft
docks with a target subject to the following discrete logic constraints:
thruster minimum impulse\dash bit, approach cone, and plume impingement. The
algorithm is demonstrated for a realistic setup inspired by the Apollo
Transposition and Docking maneuver. Fuel\dash optimal trajectories are
generated with favorable runtimes that indicate a potential for real\dash time
performance. The associated open\dash source implementation can be used as a
ground\dash based analysis tool, and can be further optimized for onboard
autonomous docking guidance.



\section*{Acknowledgments}

The authors would like to extend their gratitude to Michael Szmuk and Taylor
P. Reynolds for sharing their expertise in sequential convex programming and
for their work on the initial conference paper \cite{MalyutaScitechDocking}.

\bibliography{references}

\begin{thebibliography}{61}
\newcommand{\enquote}[1]{``#1''}
\providecommand{\natexlab}[1]{#1}
\providecommand{\url}[1]{\texttt{#1}}
\providecommand{\urlprefix}{URL }
\expandafter\ifx\csname urlstyle\endcsname\relax
  \providecommand{\doi}[1]{\discretionary{}{}{}https://doi.org/#1}\else
  \providecommand{\doi}[1]{\discretionary{}{}{}\urlstyle{rm}\url{https://doi.org/#1}}\fi

\bibitem[{Woffinden and Geller(2007)}]{Woffinden2007}
Woffinden, D.~C., and Geller, D.~K., \enquote{Navigating the road to autonomous
  orbital rendezvous,} \emph{Journal of Spacecraft and Rockets}, Vol.~44,
  No.~4, 2007, pp. 898--909.
\newblock \doi{10.2514/1.30734}.

\bibitem[{Goodman(2006)}]{Goodman2006}
Goodman, J.~L., \enquote{History of space shuttle rendezvous and proximity
  operations,} \emph{Journal of Spacecraft and Rockets}, Vol.~43, No.~5, 2006,
  pp. 944--959.
\newblock \doi{10.2514/1.19653}.

\bibitem[{D'Souza et~al.(2007)D'Souza, Hanak, Spehar, Clark, and
  Jackson}]{DSouza2007}
D'Souza, C.~N., Hanak, F.~C., Spehar, P., Clark, F.~D., and Jackson, M.,
  \enquote{Orion rendezvous, proximity operations, and docking design and
  analysis,} \emph{AIAA Guidance, Navigation and Control Conference}, Hilton
  Head, SC, 2007, pp. 1--13.
\newblock \doi{10.2514/6.2007-6683}.

\bibitem[{Weitering(2020)}]{SpaceXAutoDock}
Weitering, H., \enquote{{SpaceX}'s 1st upgraded {Dragon} cargo ship docks
  itself at space station with science, goodies and new airlock,}
  \url{https://www.space.com/spacex-cargo-dragon-crs-21-docks-at-space-station},
  Dec. 2020.

\bibitem[{Nishida and Kawamoto(2011)}]{Nishida2011}
Nishida, S.-I., and Kawamoto, S., \enquote{Strategy for capturing of a tumbling
  space debris,} \emph{Acta Astronautica}, Vol.~68, No. 1-2, 2011, pp.
  113--120.
\newblock \doi{10.1016/j.actaastro.2010.06.045}.

\bibitem[{Kaplan(2009)}]{Kaplan2009}
Kaplan, M., \enquote{Survey of Space Debris Reduction Methods,} \emph{{AIAA}
  {SPACE} 2009 Conference {\&} Exposition}, American Institute of Aeronautics
  and Astronautics, 2009, pp. 1--11.
\newblock \doi{10.2514/6.2009-6619}.

\bibitem[{Chang(2020)}]{SpaceJunkRendezvous}
Chang, K., \enquote{An Orbital Rendezvous Demonstrates a Space Junk Solution,}
  \url{https://www.nytimes.com/2020/02/26/science/mev-1-northrop-grumman-space-junk.html},
  Feb. 2020.

\bibitem[{Malyuta(2021)}]{OurOpenSourceCode}
Malyuta, D., \enquote{{SCP for Trajectory Optimization},}
  \url{https://github.com/dmalyuta/scp_traj_opt/tree/jgcd}, 2021.

\bibitem[{Malyuta et~al.(2021{\natexlab{a}})Malyuta, Reynolds, Szmuk, Lew,
  Bonalli, Pavone, and Acikmese}]{CSM2021}
Malyuta, D., Reynolds, T.~P., Szmuk, M., Lew, T., Bonalli, R., Pavone, M., and
  Acikmese, B., \enquote{Convex Optimization for Trajectory Generation,}
  \emph{{IEEE} Control Systems Magazine (work in progress)},
  2021{\natexlab{a}}.

\bibitem[{Miele et~al.(2007{\natexlab{a}})Miele, Weeks, and
  Ciarci{\`{a}}}]{Miele2007}
Miele, A., Weeks, M.~W., and Ciarci{\`{a}}, M., \enquote{Optimal trajectories
  for spacecraft rendezvous,} \emph{Journal of Optimization Theory and
  Applications}, Vol. 132, No.~3, 2007{\natexlab{a}}, pp. 353--376.
\newblock \doi{10.1007/s10957-007-9166-4}.

\bibitem[{Miele et~al.(2007{\natexlab{b}})Miele, Ciarci{\`{a}}, and
  Weeks}]{Miele2007a}
Miele, A., Ciarci{\`{a}}, M., and Weeks, M.~W., \enquote{Guidance trajectories
  for spacecraft rendezvous,} \emph{Journal of Optimization Theory and
  Applications}, Vol. 132, No.~3, 2007{\natexlab{b}}, pp. 377--400.
\newblock \doi{10.1007/s10957-007-9165-5}.

\bibitem[{Pascucci et~al.(2017)Pascucci, Szmuk, and \Acikmese}]{Pascucci2017}
Pascucci, C.~A., Szmuk, M., and \Acikmese, B., \enquote{Optimal real-time force
  rendering for on-orbit structures assembly,} \emph{10th International {ESA}
  Conference on Guidance, Navigation \& Control Systems}, {ESA}, 2017, pp.
  1--9.

\bibitem[{Harris and A{\c{c}}{\i}kme{\c{s}}e(2014)}]{Harris2014b}
Harris, M.~W., and A{\c{c}}{\i}kme{\c{s}}e, B., \enquote{Minimum Time
  Rendezvous of Multiple Spacecraft Using Differential Drag,} \emph{Journal of
  Guidance, Control, and Dynamics}, Vol.~37, No.~2, 2014, pp. 365--373.
\newblock \doi{10.2514/1.61505}.

\bibitem[{Lu and Liu(2013)}]{Lu2013}
Lu, P., and Liu, X., \enquote{Autonomous Trajectory Planning for Rendezvous and
  Proximity Operations by Conic Optimization,} \emph{Journal of Guidance,
  Control, and Dynamics}, Vol.~36, No.~2, 2013, pp. 375--389.
\newblock \doi{10.2514/1.58436}.

\bibitem[{Liu and Lu(2013)}]{Liu2013}
Liu, X., and Lu, P., \enquote{Robust Trajectory Optimization for Highly
  Constrained Rendezvous and Proximity Operations,} \emph{{AIAA} Guidance,
  Navigation, and Control ({GNC}) Conference}, American Institute of
  Aeronautics and Astronautics, 2013, pp. 1--18.
\newblock \doi{10.2514/6.2013-4720}.

\bibitem[{Breger and How(2008)}]{Breger2008}
Breger, L.~S., and How, J.~P., \enquote{Safe trajectories for autonomous
  rendezvous of spacecraft,} \emph{Journal of Guidance, Control, and Dynamics},
  Vol.~31, No.~5, 2008, pp. 1478--1489.
\newblock \doi{10.2514/1.29590}.

\bibitem[{Richards et~al.(2002)Richards, Schouwenaars, How, and
  Feron}]{Richards2002}
Richards, A., Schouwenaars, T., How, J.~P., and Feron, E., \enquote{Spacecraft
  trajectory planning with avoidance constraints using mixed-integer linear
  programming,} \emph{Journal of Guidance, Control, and Dynamics}, Vol.~25,
  No.~4, 2002, pp. 755--764.
\newblock \doi{10.2514/2.4943}.

\bibitem[{Sun et~al.(2019)Sun, Dai, and Lu}]{Sun2019}
Sun, C., Dai, R., and Lu, P., \enquote{Multi-phase spacecraft mission
  optimization by quadratically constrained quadratic programming,} \emph{AIAA
  Scitech Forum}, San Diego, CA, 2019, pp. 1--15.
\newblock \doi{10.2514/6.2019-1667}.

\bibitem[{Phillips et~al.(2003)Phillips, Kavraki, and
  Bedrossian}]{Phillips2003}
Phillips, J., Kavraki, L., and Bedrossian, N., \enquote{Spacecraft Rendezvous
  and Docking with Real-Time, Randomized Optimization,} \emph{{AIAA} Guidance,
  Navigation, and Control Conference and Exhibit}, American Institute of
  Aeronautics and Astronautics, 2003, pp. 1--11.
\newblock \doi{10.2514/6.2003-5511}.

\bibitem[{Hartley et~al.(2013)Hartley, Gallieri, and Maciejowski}]{Hartley2013}
Hartley, E.~N., Gallieri, M., and Maciejowski, J.~M., \enquote{Terminal
  spacecraft rendezvous and capture with {LASSO} model predictive control,}
  \emph{International Journal of Control}, Vol.~86, No.~11, 2013, pp.
  2104--2113.
\newblock \doi{10.1080/00207179.2013.789608}.

\bibitem[{Bengea and DeCarlo(2005)}]{Bengea2005}
Bengea, S.~C., and DeCarlo, R.~A., \enquote{Optimal control of switching
  systems,} \emph{Automatica}, Vol.~41, No.~1, 2005, pp. 11--27.
\newblock \doi{10.1016/j.automatica.2004.08.003},
  \urlprefix\url{https://doi.org/10.1016/j.automatica.2004.08.003}.

\bibitem[{Saranathan and Grant(2018)}]{saranathan2018relaxed}
Saranathan, H., and Grant, M.~J., \enquote{Relaxed Autonomously Switched Hybrid
  System Approach to Indirect Multiphase Aerospace Trajectory Optimization,}
  \emph{Journal of Spacecraft and Rockets}, Vol.~55, No.~3, 2018, pp. 611--621.
\newblock \doi{10.2514/1.a34012}.

\bibitem[{Taheri et~al.(2020)Taheri, Junkins, Kolmanovsky, and
  Girard}]{taheri2020novel}
Taheri, E., Junkins, J.~L., Kolmanovsky, I., and Girard, A., \enquote{A novel
  approach for optimal trajectory design with multiple operation modes of
  propulsion system, part 1,} \emph{Acta Astronautica}, Vol. 172, 2020, pp.
  151--165.
\newblock \doi{10.1016/j.actaastro.2020.02.042}.

\bibitem[{Arya et~al.(2021)Arya, Taheri, and Junkins}]{arya2021composite}
Arya, V., Taheri, E., and Junkins, J.~L., \enquote{A composite framework for
  co-optimization of spacecraft trajectory and propulsion system,} \emph{Acta
  Astronautica}, Vol. 178, 2021, pp. 773--782.
\newblock \doi{10.1016/j.actaastro.2020.10.007}.

\bibitem[{Szmuk et~al.(2020)Szmuk, Reynolds, and
  A{\c{c}}{\i}kme{\c{s}}e}]{Szmuk2020}
Szmuk, M., Reynolds, T.~P., and A{\c{c}}{\i}kme{\c{s}}e, B.,
  \enquote{Successive Convexification for Real-Time Six-Degree-of-Freedom
  Powered Descent Guidance with State-Triggered Constraints,} \emph{Journal of
  Guidance, Control, and Dynamics}, Vol.~43, No.~8, 2020, pp. 1399--1413.
\newblock \doi{10.2514/1.g004549}.

\bibitem[{Malyuta et~al.(2020)Malyuta, Reynolds, Szmuk, Acikmese, and
  Mesbahi}]{MalyutaScitechDocking}
Malyuta, D., Reynolds, T.~P., Szmuk, M., Acikmese, B., and Mesbahi, M.,
  \enquote{Fast Trajectory Optimization via Successive Convexification for
  Spacecraft Rendezvous with Integer Constraints,} \emph{{AIAA} Scitech 2020
  Forum}, American Institute of Aeronautics and Astronautics, 2020, pp. 1--24.
\newblock \doi{10.2514/6.2020-0616}.

\bibitem[{Szmuk(2019)}]{SzmukThesis}
Szmuk, M., \enquote{Successive Convexification \& High Performance Feedback
  Control for Agile Flight,} Ph.D. thesis, University of Washington, Seattle,
  WA, 2019.

\bibitem[{Reynolds et~al.(2020{\natexlab{a}})Reynolds, Szmuk, Malyuta, Mesbahi,
  A{\c{c}}{\i}kme{\c{s}}e, and Carson}]{Reynolds2020}
Reynolds, T.~P., Szmuk, M., Malyuta, D., Mesbahi, M., A{\c{c}}{\i}kme{\c{s}}e,
  B., and Carson, J.~M., \enquote{Dual Quaternion-Based Powered Descent
  Guidance with State-Triggered Constraints,} \emph{Journal of Guidance,
  Control, and Dynamics}, Vol.~43, No.~9, 2020{\natexlab{a}}, pp. 1584--1599.
\newblock \doi{10.2514/1.g004536}.

\bibitem[{Reynolds(2021)}]{ReynoldsThesis}
Reynolds, T.~P., \enquote{Computation Guidance and Control for Aerospace
  Systems,} Ph.D. thesis, University of Washington, Seattle, WA, 2021.

\bibitem[{Szmuk et~al.(2019{\natexlab{a}})Szmuk, Reynolds, Acikmese, Mesbahi,
  and Carson}]{Szmuk2019a}
Szmuk, M., Reynolds, T.~P., Acikmese, B., Mesbahi, M., and Carson, J.~M.,
  \enquote{Successive Convexification for 6-{DoF} Powered Descent Guidance with
  Compound State-Triggered Constraints,} \emph{{AIAA} Scitech 2019 Forum},
  American Institute of Aeronautics and Astronautics, 2019{\natexlab{a}}, pp.
  1--16.
\newblock \doi{10.2514/6.2019-0926}.

\bibitem[{Szmuk et~al.(2019{\natexlab{b}})Szmuk, Malyuta, Reynolds, Mceowen,
  and Acikmese}]{Szmuk2019b}
Szmuk, M., Malyuta, D., Reynolds, T.~P., Mceowen, M.~S., and Acikmese, B.,
  \enquote{Real-Time Quad-Rotor Path Planning Using Convex Optimization and
  Compound State-Triggered Constraints,} \emph{2019 {IEEE}/{RSJ} International
  Conference on Intelligent Robots and Systems ({IROS})}, {IEEE},
  2019{\natexlab{b}}, pp. 7666--7673.
\newblock \doi{10.1109/iros40897.2019.8967706}.

\bibitem[{Betts(1998)}]{Betts1998}
Betts, J.~T., \enquote{Survey of numerical methods for trajectory
  optimization,} \emph{Journal of Guidance, Control, and Dynamics}, Vol.~21,
  No.~2, 1998, pp. 193--207.
\newblock \doi{10.2514/2.4231}.

\bibitem[{Pontryagin et~al.(1986)Pontryagin, Boltyanskii, Gamkrelidze, and
  Mishchenko}]{PontryaginBook}
Pontryagin, L.~S., Boltyanskii, V.~G., Gamkrelidze, R.~V., and Mishchenko,
  E.~F., \emph{{The Mathematical Theory of Optimal Processes}}, Gordon and
  Breach Science Publishers, Montreux, 1986.
\newblock \doi{10.1201/9780203749319}.

\bibitem[{Berkovitz(1974)}]{BerkovitzBook}
Berkovitz, L.~D., \emph{Optimal Control Theory}, Springer New York, 1974.
\newblock \doi{10.1007/978-1-4757-6097-2}.

\bibitem[{NAS(1969{\natexlab{a}})}]{NASA_apollo_mass}
\emph{{CSM}/{LM} Spacecraft Operation Data Book, Volume 3: Mass Properties},
  National Aeronautics and Space Administration, {SNA-8-D-027(III) REV 2} ed.,
  1969{\natexlab{a}}.

\bibitem[{Curtis(2014)}]{CurtisBook}
Curtis, H.~D., \emph{{Orbital Mechanics for Engineering Students}},
  3\textsuperscript{rd} ed., Butterworth-Heinemann, Waltham, MA, 2014.

\bibitem[{Sol{\`a}(2017)}]{Sola2017}
Sol{\`a}, J., \enquote{Quaternion kinematics for the error-state Kalman
  filter,} \emph{CoRR}, 2017.

\bibitem[{NAS(1970)}]{NASA_csm_aoh}
\emph{{CSM}/{LM} Spacecraft Operation Data Book, Volume 1: {CSM} Data Book,
  Part 1: Constraints and Performance}, National Aeronautics and Space
  Administration, {SNA-8-D-027(I) REV 3} ed., 1970.

\bibitem[{NAS(1969{\natexlab{b}})}]{NASA_apollo_news}
\emph{Apollo {CSM} and {LM} News Reference: Reaction Control Subsystem}, Space
  Division of North American Rockwell Corp., 1969{\natexlab{b}}.

\bibitem[{NAS(1969{\natexlab{c}})}]{NASA_aoh_vol1}
\emph{Apollo Operations Handbook, Block II Spacecraft, Volume 1: Spacecraft
  Description}, National Aeronautics and Space Administration, {SM2A-03-Block
  II-(1)} ed., 1969{\natexlab{c}}.

\bibitem[{Boyd and Vandenberghe(2004)}]{BoydConvexBook}
Boyd, S., and Vandenberghe, L., \emph{Convex Optimization}, Cambridge
  University Press, Cambridge, UK, 2004.

\bibitem[{Nocedal and Wright(1999)}]{NocedalBook}
Nocedal, J., and Wright, S., \emph{Numerical Optimization}, Springer New York,
  1999.
\newblock \doi{10.1007/978-0-387-40065-5}.

\bibitem[{Achterberg and Wunderling(2013)}]{achterberg2013progress}
Achterberg, T., and Wunderling, R., \enquote{Mixed Integer Programming:
  Analyzing 12 Years of Progress,} \emph{Facets of Combinatorial Optimization},
  Springer Berlin Heidelberg, 2013, pp. 449--481.
\newblock \doi{10.1007/978-3-642-38189-8_18}.

\bibitem[{Schouwenaars et~al.(2001)Schouwenaars, Moor, Feron, and
  How}]{Schouwenaars2001}
Schouwenaars, T., Moor, B.~D., Feron, E., and How, J., \enquote{Mixed integer
  programming for multi-vehicle path planning,} \emph{2001 European Control
  Conference ({ECC})}, {IEEE}, 2001, pp. 2603--2608.
\newblock \doi{10.23919/ecc.2001.7076321}.

\bibitem[{Schouwenaars(2006)}]{Schouwenaars2006}
Schouwenaars, T., \enquote{{Safe trajectory planning of autonomous vehicles},}
  {Dissertation (Ph.D.)}, Massachusetts Institute of Technology, 2006.

\bibitem[{Malyuta and A\c{c}{\i}kme\c{s}e(2020)}]{MalyutaLCSS}
Malyuta, D., and A\c{c}{\i}kme\c{s}e, B., \enquote{Approximate Multiparametric
  Mixed-Integer Convex Programming,} \emph{{IEEE} Control Systems Letters},
  Vol.~4, No.~1, 2020, p. arXiv:1902.10994.
\newblock \doi{10.1109/lcsys.2019.2922639}.

\bibitem[{Harris(2021)}]{Harris2021}
Harris, M.~W., \enquote{Optimal Control on Disconnected Sets using Extreme
  Point Relaxations and Normality Approximations,} \emph{{IEEE} Transactions on
  Automatic Control}, 2021, pp. 1--1.
\newblock \doi{10.1109/tac.2021.3059682}.

\bibitem[{Malyuta and A{\c{c}}ikme{\c{s}}e(2020)}]{MalyutaIFAC}
Malyuta, D., and A{\c{c}}ikme{\c{s}}e, B., \enquote{Lossless Convexification of
  Optimal Control Problems with Semi-continuous Inputs,}
  \emph{{IFAC}-{PapersOnLine}}, Vol.~53, No.~2, 2020, pp. 6843--6850.
\newblock \doi{10.1016/j.ifacol.2020.12.341}.

\bibitem[{Blackmore et~al.(2012)Blackmore, \Acikmese{}, and {Carson
  III}}]{Blackmore2012}
Blackmore, L., \Acikmese{}, B., and {Carson III}, J.~M., \enquote{Lossless
  convexification of control constraints for a class of nonlinear optimal
  control problems,} \emph{Systems {\&} Control Letters}, Vol.~61, No.~8, 2012,
  pp. 863--870.
\newblock \doi{10.1016/j.sysconle.2012.04.010}.

\bibitem[{Harris(2014)}]{HarrisThesis}
Harris, M.~W., \enquote{Lossless Convexification of Optimal Control Problems,}
  Ph.D. thesis, The University of Texas at Austin, Austin, 2014.

\bibitem[{Malyuta et~al.(2021{\natexlab{b}})Malyuta, Yu, Elango, and
  \Acikmese{}}]{ARC2021}
Malyuta, D., Yu, Y., Elango, P., and \Acikmese{}, B., \enquote{Advances in
  Trajectory Optimization for Space Vehicle Control,} \emph{Annual Reviews in
  Control}, 2021{\natexlab{b}}.
\newblock Accepted.

\bibitem[{Reynolds et~al.(2020{\natexlab{b}})Reynolds, Malyuta, Mesbahi,
  A\c{c}{\i}kme\c{s}e, and {Carson III}}]{Reynolds2020Real}
Reynolds, T.~P., Malyuta, D., Mesbahi, M., A\c{c}{\i}kme\c{s}e, B., and {Carson
  III}, J.~M., \enquote{A Real-Time Algorithm for Non-Convex Powered Descent
  Guidance,} \emph{AIAA SciTech Forum}, AIAA, 2020{\natexlab{b}}, pp. 1--24.
\newblock \doi{10.2514/6.2020-0844}.

\bibitem[{Hastie et~al.(2009)Hastie, Tibshirani, and Friedman}]{Hastie2009}
Hastie, T., Tibshirani, R., and Friedman, J., \emph{The Elements of Statistical
  Learning}, 2\textsuperscript{nd} ed., Springer New York, 2009.
\newblock \doi{10.1007/978-0-387-84858-7}.

\bibitem[{Betts(2020)}]{BettsBook}
Betts, J.~T., \emph{Practical Methods for Optimal Control Using Nonlinear
  Programming}, {SIAM}, 2020.

\bibitem[{Watson(1986)}]{Watson1986}
Watson, L.~T., \enquote{Numerical Linear Algebra Aspects of Globally Convergent
  Homotopy Methods,} \emph{{SIAM} Review}, Vol.~28, No.~4, 1986, pp. 529--545.
\newblock \doi{10.1137/1028157}.

\bibitem[{Conn et~al.(2000)Conn, Gould, and Toint}]{ConnTrustRegionBook}
Conn, A.~R., Gould, N. I.~M., and Toint, P.~L., \emph{Trust Region Methods},
  SIAM, Philadelphia, PA, 2000.
\newblock \doi{10.1137/1.9780898719857}.

\bibitem[{Kochenderfer and Wheeler(2019)}]{Kochenderfer2019}
Kochenderfer, M.~J., and Wheeler, T.~A., \emph{Algorithms for Optimization},
  The MIT Press, Cambridge, Massachusetts, 2019.

\bibitem[{Bezanson et~al.(2017)Bezanson, Edelman, Karpinski, and
  Shah}]{bezanson2017julia}
Bezanson, J., Edelman, A., Karpinski, S., and Shah, V.~B., \enquote{Julia: A
  fresh approach to numerical computing,} \emph{SIAM review}, Vol.~59, No.~1,
  2017, pp. 65--98.
\newblock \doi{10.1137/141000671}.

\bibitem[{Domahidi et~al.(2013)Domahidi, Chu, and Boyd}]{domahidi2013ecos}
Domahidi, A., Chu, E., and Boyd, S., \enquote{{ECOS}: An {SOCP} solver for
  embedded systems,} \emph{2013 European Control Conference ({ECC})}, {IEEE},
  2013, pp. 3071--3076.
\newblock \doi{10.23919/ecc.2013.6669541}.

\bibitem[{Dueri et~al.(2014)Dueri, Zhang, and A{\c{c}}ikmese}]{Dueri2014}
Dueri, D., Zhang, J., and A{\c{c}}ikmese, B., \enquote{Automated Custom Code
  Generation for Embedded, Real-time Second Order Cone Programming,}
  \emph{{IFAC} Proceedings Volumes}, Vol.~47, No.~3, 2014, pp. 1605--1612.
\newblock \doi{10.3182/20140824-6-za-1003.02736}.

\bibitem[{Dueri et~al.(2017)Dueri, A{\c{c}}{\i}kme{\c{s}}e, Scharf, and
  Harris}]{Dueri2017}
Dueri, D., A{\c{c}}{\i}kme{\c{s}}e, B., Scharf, D.~P., and Harris, M.~W.,
  \enquote{Customized Real-Time Interior-Point Methods for Onboard
  Powered-Descent Guidance,} \emph{Journal of Guidance, Control, and Dynamics},
  Vol.~40, No.~2, 2017, pp. 197--212.
\newblock \doi{10.2514/1.g001480}.

\end{thebibliography}

\end{document}